\magnification=1200
\input amstex
\documentstyle{amsppt}

\vcorrection{-.4in}\hcorrection{.2 in}

\input xypic
\input epsf

\define\CDdashright#1#2{&\,\mathop{\dashrightarrow}\limits^{#1}_{#2}\,&}
\define\CDdashleft#1#2{&\,\mathop{\dashleftarrow}\limits^{#1}_{#2}\,&}

\def\P{\Bbb P}

\def\A{\Bbb A}

\def\G{\Bbb G}
\def\loc#1{{#1}^{\text{loc}}}

\def\Til#1{{\widetilde{#1}}}
\def\codim{\text{codim}}

\def\PGL{\text{PGL}}

\def\Spec{\text{Spec}}

\def\siltable#1.{
\vbox{\tabskip=0pt \offinterlineskip
\halign to360pt{\strut##& ##\tabskip=1em plus2em&
  \hfil##\hfil& \vrule##&
  \hfil##\hfil& \vrule##&
  \hfil##\hfil& \vrule##&
  \hfil##\hfil& \vrule\thinspace\vrule##&
  \hfil##\hfil& \vrule\thinspace\vrule##&
  \hfil##\hfil& ##\tabskip=0pt\cr
#1}}}

\CompileMatrices

\topmatter
\title Plane curves with small linear orbits I\endtitle
\author Paolo Aluffi${}^1$ and Carel Faber\endauthor
\date May 1999\enddate
\address Mathematics Department, Florida State University,
Tallahassee, FL 32306 \endaddress
\email aluffi\@math.fsu.edu \endemail
\address Department of Mathematics, Oklahoma State University,
Stillwater, OK 74078\endaddress
\email cffaber\@littlewood.math.okstate.edu \endemail

\abstract The `linear orbit' of a plane curve of degree $d$ is its
orbit in $\P^{d(d+3)/2}$ under the natural action of $\PGL(3)$. In
this paper we compute the degree of the closure of the linear orbits
of most curves with positive dimensional stabilizers. Our tool is a
nonsingular variety dominating the orbit closure, which we construct
by a blow@-up sequence mirroring the sequence yielding an embedded
resolution of the curve.

The results given here will serve as an ingredient in the computation
of the analogous information for arbitrary plane curves. Linear
orbits of smooth plane curves are studied in \cite{A-F1}.
\endabstract


\leftheadtext{Paolo Aluffi and Carel Faber}

\endtopmatter

\document
\footnote[]{${}^1$Supported in part by NSF grant DMS-9500843}

\head \S0. Introduction\endhead

In this paper we study the `linear orbits' of certain singular plane
curves. We have dealt with orbits of {\it smooth\/} plane curves in
\cite{A-F1}; the results in this paper are the next natural step
towards a treatment of {\it arbitrary\/} plane curves.

Here is the set-up. The group $\PGL(3)$ of projective transformations
of the plane $\P^2$ acts naturally on the projective space $\P^N$
parametrizing plane curves of degree $d$ (here $N=\frac{d(d+3)}2$).
The orbit of a curve $C$ is a quasi-projective variety of
dimension$\le 8$, which we call the `linear orbit' of $C$. Most curves
have linear orbits of dimension 8; we say that $C$ has a {\it small\/}
linear orbit if the dimension of its orbit is 7 or less. This paper
studies the enumerative geometry of most plane curves whose orbit is
small.

It is natural to study the closures of these linear orbits in the
projective space $\P^N$: questions
arise as to e.g.~the degrees of these projective varieties (on what
features of a plane curve does the degree of its orbit closure
depend?); the decomposition of their boundaries in smaller orbits;
their singularities (which orbit closures are smooth?); and the
behavior of orbit closures in families of plane curves.

In \cite{A-F1} we answer some of these questions in the case of a
smooth plane curve. Our main tool is the construction through explicit
blow@-ups of a nonsingular projective variety dominating the orbit
closure. The degree of the orbit closure can then be determined with
the aid of standard intersection theory. The answer depends naturally
on the degree of the plane curve and the order of its stabilizer, but
also (somewhat surprisingly) on the types of its flexes: in fact, the
structure of the blow@-up sequence depends precisely on the number and
type of the flexes on the curve.

Unfortunately, this natural approach seems inadequate for most
singular curves: we do not know a sequence of blow@-ups producing a
nonsingular variety dominating the orbit closure for arbitrary
singularities. In a different approach that we have developed for the study
of orbit closures, the first step is to determine which orbits appear
in the boundary of the orbit closure of a given curve; this was in
essence carried out more than 60 years ago in \cite{Ghizzetti}
\footnote{We are grateful to the referee of \cite{A-F1} for pointing
us to Ghizzetti's work.}, 
and will be discussed elsewhere.
The second step is
to study these `small' orbits in detail; the present paper contains
such a study, for almost all small orbits. More precisely, we deal
here with all curves whose orbit is small and which contain some
non@-linear component. Curves consisting entirely of lines require a
different (and in some sense simpler) treatment; their orbits, and the
classification of small orbits, are the subject matter of \cite{A-F4}.

In \S1 we describe the curves that we study in this paper, and state
the main result: the computation of the degree of the orbit closures
of these curves. These degrees (together with the related results
of \S4) will be the input necessary to treat
arbitrary singular curves. 

For curves with small orbits, the precise knowledge of the
singularities that can arise allows us to carry through the approach
used for smooth curves. The computation is again based upon the
construction (\S2) of a non-singular projective variety admitting a
dominant morphism to the orbit closure. The explicit blow-up sequence
yielding this variety now mimics the embedded resolution of the
singular curve in the plane (as mentioned above, this approach
surprisingly does {\it not\/} seem to work for arbitrary singular curves).

In \S3 we describe the actual degree computation, which is rather
involved; the main tool is a refinement (Proposition~2.3) of a
blow@-up formula from \cite{Aluffi}. The final answer (Theorem~1.1) has a
remarkably simple form, considering the laborious
procedure leading to it. For example, while the blow@-up sequence
we use relies in an essential way on the Dynkin diagrams of the
singularities,
only very coarse numerical information (such as the degree of the
components of the curve, or their multiplicity at the singular points)
enters in the formula for the degree of its orbit closure.

In \S4 we also discuss `predegree polynomials', which combine
information concerning the enumerative geometry of the curves when
certain natural constraints are 
introduced. Formulas for the degrees of loci of curves with these
constraints are obtained by applying a suitable differential operator
to the expression in Theorem~1.1.
These results are included both because they are natural extensions of
the other results in this paper, and because they will be ingredients
in the computation of the degree of the orbit closure of an {\it
arbitrary\/} plane curve, which we will describe elsewhere.
\vskip 6pt

{\bf Acknowledgements.\/} We thank the University of Chicago, 
Mathematisches Forschungsinstitut Oberwolfach, and the Mittag--Leffler
Institut, for hospitality and support; and W.~Fulton, M.~Kreck, 
and D.~Laksov for
the invitations to visit these institutions.
Our research at Oberwolfach was supported by the R.i.P.~program,
generously funded by the Volkswagen@-Stiftung.
Finally, we thank the referee for useful comments on an earlier
version of this paper.

\head \S1. Statement of the main result\endhead

We work over an algebraically closed field of characteristic 0.

Let $m<n$ be coprime integers, with $m\ge 1$. The prototype
irreducible curve we consider in this paper is the cuspidal plane
curve $C$ `of type $(m,n)$', i.e., with projective equation
$$x^n=y^m z^{n-m}$$
for suitable coordinates $(x:y:z)$. We aim to studying the locus of
all curves of type $(m,n)$, which form the $\PGL(3)$@-orbit of a single
such curve. In fact, we are interested in studying {\it all\/} curves
whose $\PGL(3)$@-orbit has dimension $<8$; so we will study here the
orbit of a more general (possibly reducible) type of curve, specified
below.

Note that type $(m,n)$ and type $(n-m,n)$ only differ by a coordinate
switch $y\leftrightarrow z$. The only two (possibly) singular points
of $C$ are located at $(0:0:1)$ and $(0:1:0)$; we will generally call
these points `cusps', although they may in fact be nonsingular (for
$m=1$ or $m=n-1$, respectively). 
Note also that $C$ determines a triangle, formed by the line
$\lambda=\{x=0\}$, joining the two cusps, and by the tangent cones
$\mu=\{y=0\}$, $\overline\mu=\{z=0\}$ to $C$ at the cusps:
$$\epsffile{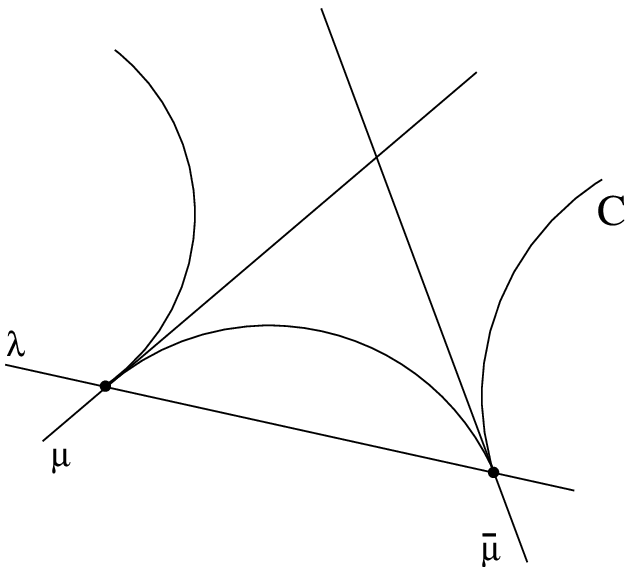}$$

More generally, fix two coprime integers $n>m\ge 1$. The curves $C$ we
study in this paper consist of arbitrary unions of curves from the pencil 
$$x^n=\alpha y^m z^{n-m}\quad(\alpha\ne 0),$$
counted with arbitrary multiplicities $s_i$, and of the lines
$\lambda$, $\mu$, $\overline \mu$ of the basic triangle, taken with
multiplicities $r$, $q$, $\overline q$ respectively. We denote by $S$
the sum $\sum s_i$, and write $\overline m=n-m$ for convenience.

Now we act on the plane by the group $\PGL(3)$ of projective linear
transformations. This action induces a (right) action on the
projective space $\P^N$, $N=\frac{d(d+3)}2$, parametrizing degree@-$d$
plane curves. A curve $C$ as specified above has degree
$d=(Sn+r+q+\overline q)$, and its orbit in $\P^N$ has dimension~7 for
all but very special cases (for example, if $S=0$ then $C$ consists of
lines from the basic triangle, and the dimension of its orbit is
necessarily $\le 6$).

In case $C$ contains, besides lines, at most one curve of type $(m,n)$,
the set of all curves of the same type and with the same multiplicities
$s_i$, $r$, $q$, $\overline q$ is precisely the orbit of $C$; we study
the closure of this orbit. If $C$ contains two or more curves from
the pencil, then the set of all curves of the same type and with the
same multiplicities consists of infinitely many orbits. 
As explained in the introduction, we study the orbits of these curves
rather than the set of all of them. We will find that the infinitely
many orbits for a given set of data have essentially the same
behavior; a special choice of the curves in the pencil may give rise
to a bigger automorphism group, which affects the degree of the orbit
closure only by a multiplicative factor.

Here is the main numerical result of the paper. First, working in the ring
with $r^3=q^3={\overline q}^3=0$, expand the expression
$$\multline
n^2 m^2 \overline m^2\left(\left(S+\frac rn+\frac qm +\frac{\overline q}
{\overline m}\right)^7 + 2\left(S+\frac
rn+\frac qm \right)^7+2\left(S+\frac rn+\frac{{\overline q}}{\overline
m}\right)^7\right.\\
\left.+\left(S+\frac rn\right)^7-42 \left(S+\frac
rn\right)^5\left(\frac{q^2}{m^2}- \frac qm \frac{{\overline q}} {\overline m}
+ \frac{{\overline q}^2}{\overline m^2}\right)^{\vphantom{7}}\right)
\endmultline$$
obtaining a {\it polynomial\/} in all the variables (of degree$\le 2$
in $r$, $q$, $\overline q$); 

\noindent---then, subtract
$$\left( 84\,{(Sn+r+q+\overline q)^2}\,{\sum s_i^5} -
252\,(Sn+r+q+\overline q) \,{\sum s_i^6} + 192\,{\sum s_i^7}
\right)\quad.$$

The result is a polynomial expression $Q(n,m,s_i,r,q,\overline q)$.

\proclaim{Theorem 1.1} If 7@-dimensional, the orbit closure of a curve
with data $n,m,s_i$, $r, q, \overline q$ as above has degree
$$\frac 1A\cdot  Q(n,m,s_i,r,q,\overline q)\quad,$$
where $A$ is the number of components of the $\PGL(3)$@-stabilizer of
the curve. If the orbit has dimension lower than 7, the expression
evaluates to 0.\endproclaim
The number $A$ accounts for special automorphisms of the curve, due to
extra symmetries in the position of the components in the pencil; 
cf.~Lemma 3.1. $A$ equals 1 for most choices of $n, m$, etc. 

The enumerative meaning of the formula obtained in Theorem~1.1
rests on the fact that imposing $C$ to contain a given point is a
linear condition in $\P^N$. If all multiplicities are 0 or 1,
it is easy to see that the number computed in the theorem equals the
number of curves in the orbit which contain~7 general points.

The first of the two expressions building up to $Q$ will be obtained
by combining a `B\'ezout term' with contributions arising from the
`local' part of the construction in \S2, essentially aimed at
resolving the singularities of the curve.
Note that it only depends on the multiplicities $s_i$ of the cuspidal
components via their sum $S$. The second term will arise from the
`global' stage of the construction, taking care of the curve after
singularities have been resolved. The multiplicities enter here in a
more interesting way, but note that this term depends otherwise only
on the total degree $d$ of the curve. We do not have
conceptual explanations for these features, or for the remarkable
shape of the first expression (indeed, our construction only yields a
complicated raw expression, which we then recognize to equal the
relatively simple one given above).

To our knowledge, there is minimal overlap of the results in this
note with the existing literature in enumerative geometry. J.~M.~Miret
and S.~Xamb\'o have computed hundreds of characteristic numbers for
cuspidal plane {\it cubics,\/} in \cite{M-X}; our formulas allow us to
reproduce 27 of the numbers in their lists. In fact, a particular case
of our result yields closed formulas for these numbers for curves of
{\it arbitrary\/} degree, in terms of the type $(m,n)$ of the curve
(see \S4.3).

Notice that the dual of a curve of type $(m,n)$
as above is a curve of type $(\overline m,n)$, hence again of type
$(m,n)$ after a coordinate switch. Therefore, the degrees computed
here also compute {\it characteristic numbers:\/} that is, the number
of curves of given type and tangent to 7~lines in general position in
the plane.

The identity component of the stabilizer of a curve with 7@-dimensional
linear orbit is either $\G_m$ or $\G_a$.
All curves containing some non-linear component and whose stabilizer
contains $\G_m$ are of the kind considered above (\cite{A-F4}).
The curves with 7@-dimensional orbit and whose stabilizer
contains a $\G_a$ are not of this kind
(one example of such curves is the union of two smooth conics
touching at exactly one point). This case is briefly discussed in
\S4.1;
the formula given above turns out to be correct for this case as well,
with suitable choices of the variables.

We also include here (see \S4.2) a few remarks that extend the result
given above to cases in which the orbit has dimension $< 7$. Moreover,
we discuss the degree of subsets of the orbit closures determined by
imposing conditions on the lines of the basic triangle (see \S4.3).

\head \S2. Local and global blow@-ups\endhead

Our goal in this section is the explicit construction of nonsingular
varieties dominating the closure $\overline{\Cal O}_C$ of the orbits
of the curves discussed in \S1. The general approach we take is
a natural extension of the one in \cite{A-F1}, and we summarize it here.

After choosing coordinates in $\P^2$, we consider the $\P^8$ of
$3\times 3$ matrices as a completion of $\PGL(3)$. The action of
$\PGL(3)$ on a fixed curve $C$ determines then a rational map
$$\P^8 \overset c \to\dashrightarrow \P^N\quad,$$
by sending a matrix $\varphi$ to the curve with equation
$F(\varphi(x:y:z))=0$, where $F$ is an equation for $C$. Our aim is to
resolve the
indeterminacies of this map, by a sequence of blow@-ups at nonsingular
centers, starting from $\P^8$. We will then obtain a {\it
nonsingular\/}
variety $\Til V$ surjecting onto the orbit closure:
$$\Til V @>>> \overline{\Cal O}_C\quad.$$
The challenge is to perform the resolution explicitly enough to be
able to keep track of the intersection theory and of other relevant
information.
If this is
accomplished, then several invariants of $\overline{\Cal O}_C$ (such
as degree, Euler characteristic, multiplicity along components of the
singular locus, etc.) can be computed in principle. This will be
illustrated in \S3 by the computation of the degree of $\overline{\Cal
O}_C$, with the result stated in \S1. The computation of other
invariants might be substantially more involved; for an example in
which multiplicity computations can be carried out explicitly, see
\cite{A-F2}. 
With this broader range of problems in mind, we insist on aiming to
construct a {\it nonsingular\/} $\Til V$, although this forces us into
a bit of extra work in this \S 2.

As is immediately checked, the base locus of the rational map $c$
defined above is supported on the set of rank@-1 matrices whose image
is a point of $C$, union the set of rank@-2 matrices whose image is a
line contained in $C$.

The resolution of the indeterminacies of $c$ will require two distinct
stages. In a first stage we will deal with the fact that the curves we
consider are {\it singular\/} (in general): this causes the base locus
of $c$ to be itself singular, and we employ a sequence of blow@-ups to
resolve its singularity. We call these blow@-ups `local', to remind
ourselves that they deal with local features of the curves under exam.
Once the singularities of the base locus are resolved, we need a
second stage of `global' blow@-ups to eliminate the indeterminacies of
the lifted rational map.  This stage is considerably simpler,
particularly because the situation is reduced to the case
of nonsingular curves, which was examined in \cite{A-F1}.

The details of the construction are rather technical; however, a
rather explicit description of a variety $\Til V$ as above is
necessary in order to perform the degree computations in \S3, and
would be essential to attack subtler problems such as the
study of singularities of the orbit closure. We therefore feel that it
would not be opportune to omit these details altogether. Here is a
summary of how the section is organized; the hurried reader should
feel free to skip the rest of this \S2 at first reading.

---We consider a curve $C=$ the union of finitely many curves of type
$(m,n)$, and of lines from the basic triangle, with arbitrary
multiplicities (see \S1);

---The action of $\PGL(3)$ extends to a dominant rational map from the
$\P^8$ of $3\times 3$ matrices to the orbit closure of $C$: $c:\P^8
\dashrightarrow \overline{\Cal O}_C\subset\P^N$;

---The indeterminacies of this map are removed by a sequence of
blow@-ups, and more precisely:\roster
\item"(i)" a sequence of `local' blow@-ups of $\P^8$ along nonsingular
centers, corresponding to the two cusps of the components of $C$.
These blow@-ups mirror the sequence of blow@-ups yielding the
embedded resolution of $C$. The sequence corresponding to one cusp is
described in Proposition 2.6. This produces a variety
$\loc{\Til V}$;
\item"(ii)" two `global' blow@-ups with nonsingular centers of
dimension 3, 4, over $\loc{\Til V}$; these are discussed in Theorem 2.4;
\item"(iii)" a blow@-up along a $\P^2$ obtained as the (isomorphic)
inverse image of the set of matrices whose image is
$\mu\cap\overline\mu$; and 
\item"(iv)" blow@-ups along three 5@-dimensional nonsingular
varieties. These are the proper transforms of the set of matrices
whose image is one of the lines of the basic triangle, see the
discussion in \S2.3.
\endroster
We denote by $\Til V$ the variety obtained at the end of this process.
By pasting together the pieces of our discussion, we will have:
\proclaim{Theorem 2.1} The procedure described above produces a
nonsingular variety $\Til V$ mapping to $\P^8$, such that $c$ lifts to
a regular map $\tilde c:\Til V @>>> \P^N$. The image of the map
$\tilde c$ is the orbit closure ${\overline {\Cal O}}_C$:
$$\diagram
{\Til V} \dto \drto^{\tilde c}\\
{\P^8} \xdashed[0,1]^{c}|>\tip & {\overline{\Cal O}_C} & {\hskip
-30pt\subset \P^N} \\
\enddiagram$$
\endproclaim
This section is devoted to the construction of $\Til V$ and the proof
of Theorem 2.1.

\subheading{\S 2.1. Directed blow@-ups}
As mentioned above, the local blow@-ups will essentially mirror the
blow@-ups needed to obtain an embedded resolution of the cuspidal
curve of type $(m,n)$ presented in \S1. We start by recalling how this
resolution is accomplished, and introduce a device (`directed
blow@-ups') which streamlines the construction considerably.

Consider the affine portion of a cuspidal curve $C$ of type $(m,n)$,
centered at one of the cusps:
$$x^n=\alpha y^m\quad (\alpha\ne 0),$$
together with the tangent cone $\mu$ to $C$ at the cusp:
$$\epsffile{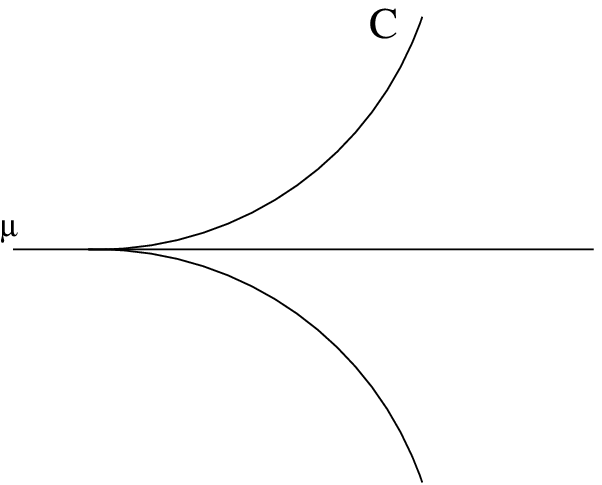}$$

To get an embedded resolution of the union $C\cup \mu$, start by
blowing up the plane at the cusp $C\cap\mu$, then successively
blow@-up at the point of intersection of the proper transform of $C$
with the latest exceptional divisor, until the resolution is achieved.
A more refined description of this sequence is controlled by the steps
of the Euclidean algorithm for $m=m_1$, $n$:
$$\align
n&=m_1 \ell_1+m_2\\
m_1&=m_2 \ell_2 +m_3\\
&\cdots\\
m_{e-2}&=m_{e-1} \ell_{e-1}+m_e\\
m_{e-1}&=m_e \ell_e
\endalign$$
with all $\ell_i$, $m_i$ positive integers, $m_i<m_{i-1}$, and $m_e=1$.

The center of each of the first $\ell_1$ blow@-ups is the
intersection of the proper transforms of $C$ and $\mu$; if $m_1=1$,
these $\ell_1$ blow@-ups produce the resolution. If $m_1>1$, after
this sequence the proper transform of $C$ is\roster
\item"(i)" transversal to the proper transform of $\mu$; and
\item"(ii)" an affine cuspidal curve of type $(m_2,m_1)$ (in a
suitable chart of the blow@-up), with tangent cone at the cusp equal
to the latest exceptional divisor $E_1$.\endroster

That is, at the end of the first $\ell_1$ blow@-ups we are left with the
same problem with which we had started, but related to a `simpler'
curve, of type $(m_2,m_1)$.

Similarly, the second line of the Euclidean algorithm corresponds to a
sequence of $\ell_2$ blow@-ups, at the end of which the proper transform
of $C$ will be a type@-$(m_3,m_2)$ curve with the latest exceptional
divisor, $E_2$, as tangent cone (if $m_2\ne 1$). Proceeding in this
fashion, the situation simplifies until the last $\ell_e$ blow@-ups,
which yield a curve `of type $(0,1)$'---that is, a curve transversal to
the last exceptional divisor $E_e$. At this point the embedded
resolution is achieved.

Note that this subdivision of the resolution process in $e$ steps,
according to the lines of the Euclidean algorithm, is natural from the
point of view of the multiplicity of the curve at the successive
centers of blow@-ups: this is $m_1$ for the $\ell_1$ blow@-ups
corresponding to the first line, then $m_2$ for the next $\ell_2$
blow@-ups, etc.

In view of these considerations, and of how they will be mirrored by
the `local' blow@-ups over $\P^8$, we define a notion of `directed'
blow@-up. Let
$$B\subset P\subset V$$
be three nonsingular varieties, with $\dim B <\dim P<\dim V$, and
let $\ell>0$. We define a {\it nonsingular\/} variety $V^{(\ell)}$
birational to $V$, and dominating the blow@-ups of $V$ along the
`$j$@-th thickening of $B$ in the direction of $P$' for all $0\le j\le
\ell$ (see the example following Lemma 2.2).

\definition{Definition} With $B\subset P\subset V$ as above, we let
$V^{(1)}$ be the blow@-up of $V$ along $B$; for $\ell\ge 2$, we let
$V^{(\ell)}$ be the blow@-up of $V^{(\ell-1)}$ along the intersection
of the proper transform of $P$ with the exceptional divisor
$E^{(\ell-1)}$ in $V^{(\ell-1)}$. We call $V^{(\ell)}$ the {\it
$\ell$@-directed blow@-up of $V$ along $B$, in the direction of
$P$.\/}
The {\it exceptional divisor\/} of the directed blow@-up is the last
exceptional divisor, $E^{(\ell)}$, produced in the sequence.
Also, the exceptional divisor of the directed blow@-up contains a
{\it distinguished\/} subvariety, namely its intersection with the
proper transform of $P$.
\enddefinition

With this terminology, each stage of the resolution described above
(corresponding to one line of the Euclidean algorithm) is simply one
directed blow@-up at the cusp, in the direction of the tangent cone.

Directed blow@-ups satisfy a few simple properties, whose proof we
leave to the reader:
\proclaim{Lemma 2.2} \roster
\item For $j\ge 1$, the proper transform  $\Til P\subset V^{(j)}$
of $P$
is isomorphic to the blow@-up $\Til P=B\ell_BP$ of $P$ along $B$. The
centers of the blow@-ups in $V^{(j)}$ are all isomorphic to the
projectivization $\Til B$ of the normal bundle $N_BP$ of $B$ in $P$. Let
$\Cal O(-1)$ denote the universal line subbundle on $\Til B=\P(N_BP)$.
\item For all $j\ge 1$, $\Til P$ and $\Til B\subset E^{(j)}$ are
disjoint from the proper transforms of `previous' exceptional
divisors $E^{(i)}$, $i<j$. Also, 
$$E^{(j)}\cdot \Til B=c_1(\Cal O(-1))\cap [\Til B]\quad.$$
\item For all $j\ge 1$, 
$$c(N_{\Til B}V^{(j)})=c(\Cal O(-1))\cdot c(N_PV\otimes \Cal O(-j))
\quad.$$
\item Let $(x_1,\dots,x_n)$ be local parameters for $V$, such that $P$
and $B$ are locally given respectively by the ideals $(x_1,\dots,x_p)$
and $(x_1,\dots,x_b)$ ($p<b$). Then near $\Til B\subset V^{(\ell)}$ we may
give covering charts $U_{p+1},\dots, U_b$ and local parameters
$(y_1,\dots,y_n)$ on $U_j$ so that the composition of the blow@-up maps
$$U_j @>>> V$$
is given by
$$x_i=\left\{\aligned
y_i y_j^\ell\quad\quad &i=1,\dots,p\\
y_i y_j \quad\quad &i=p+1,\dots,\hat j,\dots,b\\
y_i  \quad\quad &i=j, i=b+1,\dots,n
\endaligned\right.$$
The ideal of $E^{(\ell)}$ in this chart is $(y_j)$; the ideal of $\Til P$
is
$$(y_1,\dots,y_p)$$
\item Let $\Cal J\subset\Cal I$ be the ideal sheaves of $P, B$
respectively. For all $j\ge 1$, let $B^{(j)}$ be the subscheme of $V$
defined by the ideal $(\Cal I^j+\Cal J)$. Then
$B=B^{(1)},\dots,B^{(\ell)}$ all pull@-back to Cartier divisors
on $V^{(\ell)}$.
\endroster\endproclaim

\example{Example} To clarify the construction, let us compare the
directed blow@-up corresponding to the first line of the Euclidean
algorithm for the curve $x^n=y^m$ (with $n=\ell_1 m+m_2$, $B=$the
origin, and $P=$the $x$@-axis; so, $\Cal I=(x,y)$ and $\Cal J=(y)$)
with blowing@-up the plane directly along the fat point with ideal
$\Cal I^{\ell_1}+\Cal J=(x^{\ell_1},y)$.
The $\ell_1$@-directed
blow@-up produces $\ell_1$ exceptional divisors, and (as we
observed already) the proper transform of the curve is the curve
$t^m-x^{m_2}$ in a suitable chart of the resulting nonsingular surface.
The blow@-up along $\Cal I^{\ell_1}+\Cal J$ is covered by charts
$$\Spec\frac{k[x,y,t]}{(y-t x^{\ell_1})}\quad,\quad
\Spec\frac{k[x,y,s]}{(sy-x^{\ell_1})}\quad;$$
so the total transform of the curve $x^n=y^m$ is covered by
$$\align
\Spec\frac{k[x,y,t]}{(y-t x^{\ell_1},x^n-y^m)} &\cong
\Spec\frac{k[x,t]}{(x^{\ell_1 m}(t^m-x^{m_2}))}\\
\Spec\frac{k[x,y,s]}{(sy-x^{\ell_1},x^n-y^m)} &\cong
\Spec\frac{k[x,y,s]}{(sy-x^{\ell_1},y^m(1-s^m x^{m_2}))}
\endalign$$
{From} this we see that the proper transform of the curve sits in the
nonsingular part of the blow@-up, and there it behaves just as in the
$\ell_1$@-directed blow@-up of the plane along the origin, in the
direction of the tangent cone.\endexample

The $\ell$@-directed blow@-up of $V$ along $B$ in the direction of
$P$ is simply a resolution of singularities of the blow@-up along the
subscheme $\Cal I^\ell+\Cal J$. The few extra exceptional divisors
introduced in the process seem a price worth paying for the benefit of
obtaining a {\it nonsingular\/} variety dominating the orbit closure. 

The `local blow@-ups' in the resolution of the map $c$ introduced at
the beginning of this section will be a sequence of directed
blow@-ups, also controlled by the Euclidean algorithm on $m,n$. The
application in \S3 (yielding the degree of the closure of the image of
$c$) will rely on keeping track of the intersection of several
divisors in these blow@-ups. We will make use of the following
formula, which compares the intersection number of a collection of
divisors with the intersection number of their proper transforms after
a directed blow@-up.

To state this formula, we will use the notation $c(\Cal E^{(\ell)})$ for
the Adams operation on the Chern class of a bundle $\Cal E$:
$$c_r(\Cal E^{(\ell)})=\ell^r c_r(\Cal E)\quad,\quad r\ge 0\quad.$$

\proclaim{Proposition 2.3} Denote by $V^{(\ell)} @>\pi>> V$ the
$\ell$@-directed
blow@-up of $V$ along $B$ in the direction of $P$, as above.
Let $X_i$, $i=1,\dots,j$ and $Y_i$, $i=1,\dots,k$ be effective
divisors in $V$, such that

---each $X_i$ and its proper transforms contain the centers $B$, $\Til B$
of blow@-ups with the same multiplicity $m_i$, and

---each $Y_i$ has multiplicity $r_i$ along $B$, and its proper
transforms do not contain the other centers $\Til B$ of blow@-ups.

Further, assume that the number $k$ of divisors $Y_i$ be less than the
codimension of $B$ in $P$. Denote by $\Til X_i$, $\Til Y_i$ the proper
transforms of the divisors in $V^{(\ell)}$. Then
$$\prod_{i=1}^k Y_i\cdot \prod_{i=1}^j X_i \cap[V]-
\pi_*\left(\prod_{i=1}^k \Til Y_i\cdot \prod_{i=1}^j \Til X_i \cap
[V^{(\ell)}]\right)$$
equals the push@-forward from $B$ of the term of dimension $\dim
V-j-k$ in
$$ \ell^{\dim P-\dim B-k} \frac{\prod_{i=1}^k (r_i+\ell Y_i)
\prod_{i=1}^j (m_i+X_i)}{c(N_B^{(\ell)}P)c(N_PV)} \cap [B]\quad.$$
\endproclaim

For $\ell=1$, that is for the ordinary blow@-up of $V$ along $B$, this
says that the intersection of the divisors changes under proper
transforms by the term of expected dimension in
$$ \frac{\prod (r_i+Y_i) \prod (m_i+X_i)}
{c(N_BV)} \cap [B]\quad.$$
This is a restatement of a particular case of Theorem~II in
\cite{Aluffi}. The formula for directed blow@-ups can be deduced from
the $\ell=1$ case; we leave the details to the reader.

\subheading{\S2.2. Blow@-ups for one curve}
We first consider the case of a single irreducible curve $C$ of type
$(m,n)$, and focus our attention on one of the cusps. So choose affine
coordinates and write the equation of $C$
$$x^n=\alpha y^m$$
with $\alpha\ne 0$, $1\le m<n$ and $(m,n)=1$. We may and will in fact
assume $\alpha=1$, by rescaling $y$. Note that the base locus of the
corresponding rational map $c:\P^8 \dashrightarrow \P^N$ contains the
set $B\cong \P^2$ of rank@-1 matrices whose image is the cusp of $C$.
In fact, the base locus of $c$ consists of an isomorphic copy of
$\P^2\times C$: the set of rank@-1 matrices with arbitrary kernel, and
image a point on $C$. For $m>1$, this locus is singular along $B$.

As mentioned in the summary preceding the statement of Theorem 2.1, we
construct our resolution $\Til V$ by a two@-stage process. The first
stage consists of a sequence of directed blow@-ups, mirroring the
sequence giving the resolution of the union of the curve and its
tangent cone at the cusp. More precisely, assume that the Euclidean
algorithm for $m, n$ consists of $e$ lines, as in \S2.1:
$$\align
n&=m_1 \ell_1+m_2\\
m_1&=m_2 \ell_2 +m_3\\
&\cdots\\
m_{e-2}&=m_{e-1} \ell_{e-1}+m_e\\
m_{e-1}&=m_e \ell_e
\endalign$$

Also, the base locus of $c$ is a copy of $\P^2\times C$; the cuspidal
point $(0,0)$ of $C$ determines a distinguished
$B=\P^2\times\{(0,0)\}$ in the base locus. The tangent cone to
$C$ at $(0,0)$ is the line $y=0$, which determines a distinguished
$P=\P^5\subset\P^8$, that is the set of matrices whose image is
contained in this line.

\definition{Definition} We define a variety $\loc{\Til V}$ by the
following sequence of $e$ directed blow@-ups:

---first, perform the $\ell_1$@-directed blow@-up of $\P^8$ along $B$
in the direction of $P$; this produces a variety $\loc V_1$, with an
exceptional divisor $\loc E_1$ and a distinguished 4@-fold $\loc B_1
\subset \loc E_1$;

---next, for $i>1$ perform inductively the $\ell_i$@-directed
blow@-up of $\loc V_{i-1}$ along $\loc B_{i-1}$ in the direction of
$\loc E_{i-1}$; this produces a variety $\loc V_i$, with an
exceptional divisor $\loc E_i$ and a distinguished 6@-fold $\loc
B_i\subset \loc E_i$;

with these notations, we let $\loc{\Til V}=\loc V_e$.\enddefinition

In order to study $\loc{\Til V}$, and to describe the second stage of
the process, we introduce affine coordinates
$$\pmatrix 1 & p_1 & p_2 \\
p_3 & p_4 & p_5 \\
p_6 & p_7 & p_8
\endpmatrix$$
for $\P^8$ near $\pmatrix 1 & 0 & 0 \\0 & 0 & 0\\ 0 & 0 &
0\endpmatrix$, that is a rank@-one matrix with image the origin
$(1:0:0)$, and remark that this choice of coordinates is irrelevant in
the sense that we can move (by multiplying on the right by a constant
matrix) any such $3\times 3$ matrix to one in the chosen $\A^8$, so
that we must be able to detect in this $\A^8$ every phenomenon
relevant to our computation. 

The sequence of `local' directed blow@-ups specified above produces a
variety $\loc{\Til V} @>{\loc\pi}>> \P^8$; we give coordinates on a chart
in $\loc{\Til V}$:
$$\pmatrix
\boxed{\loc{\Til V}} &s_1 & s_2\\
s_3 & s_4 & s_5\\
s_6 & s_7 & s_8
\endpmatrix$$
(the boxed entry reminds us of the variety where the coordinates are
given) and we claim that again the choice of this particular chart will be
irrelevant for what follows. The expression of $\loc \pi:\loc{\Til V}
@>>> \P^8$ in these coordinates depends on the parity of the number of
steps $e$ in the Euclidean algorithm displayed above:

if $e$ is {\it odd,\/} we will have
$$\pmatrix
1 & p_1 & p_2 \\
p_3 & p_4 & p_5 \\
p_6 & p_7 & p_8
\endpmatrix=\pmatrix
1 & s_1 & s_2 \\
s_3^m s_6^A & s_3^m s_6^As_4 & s_3^m s_6^As_5 \\
s_3^n s_6^B & s_3^n s_6^Bs_7 & s_3^n s_6^Bs_8 \\
\endpmatrix\quad,$$
with $Bm-An=1$ (the actual values of $A$, $B$ can be obtained in terms
of the Euclidean algorithm, but are not important here);

if $e$ is {\it even,\/} we will have
$$\pmatrix
1 & p_1 & p_2 \\
p_3 & p_4 & p_5 \\
p_6 & p_7 & p_8
\endpmatrix=\pmatrix
1 & s_1 & s_2 \\
s_3^A s_6^m & s_3^A s_6^ms_4 & s_3^A s_6^ms_5 \\
s_3^B s_6^n & s_3^B s_6^ns_7 & s_3^B s_6^ns_8 \\
\endpmatrix\quad,$$
with $An-Bm=1$.

\remark{Remark} These coordinate expressions are slightly different in
the case $m=1$, i.e., $e=1$. Other details of the construction
require minor modifications in this case; we leave
these to the reader.\endremark

In Proposition 2.6 we will prove that coordinates can be given on
$\loc{\Til V}$ so that these expressions hold. First, 
we claim that if we show that this coordinate description holds (and
that the choice of the chart is indeed irrelevant), then we are
essentially done:

\proclaim{Theorem 2.4} Two blow@-ups at smooth centers remove the
indeterminacies of the lifted map $\loc c: \loc{\Til V}
\dashrightarrow \P^N$.\endproclaim
\proclaim{Corollary 2.5} For the curve $C$ with equation $x^n=\alpha y^m
z^{n-m}$, $\alpha\ne 0$, the indeterminacies of the corresponding
rational map $c$ can be removed by performing the sequence of `local'
blow@-ups for the two cusps, followed by two `global' blow@-ups at
smooth centers.\endproclaim

We prove the theorem right away, and concentrate on the more involved
details of the coordinate description in Proposition 2.6. By
a `point@-condition' we mean the hypersurface in $\P^8$ formed by all
matrices which send the chosen curve $C$ (say with equation $F=0$) to
contain a fixed point $p$. More precisely, the point@-condition in
$\P^8$ corresponding to $p\in \P^2$ has equation (in $\varphi\in\P^8$)
$$F(\varphi(p))=0\quad.$$
Further, we call `point@-conditions' the proper transforms of
point@-conditions in any variety mapping birationally to $\P^8$.
The point@-conditions in $\P^8$ generate the linear
system corresponding to $c$; hence, showing that a lift of
$c$ to a variety $\Til V$ removes the indeterminacies of $c$ amounts to
showing that the (proper transforms of the) point@-conditions in $\Til
V$ do not have a common intersection.

\demo{Proof of Theorem 2.4} The point@-condition corresponding to
$(\xi_0, \xi_1,\xi_2)$ in $\A^8\subset \P^8$ has equation
$$(\xi_0+p_1\xi_1+p_2\xi_2)^{n-m}(p_6\xi_0+p_7\xi_1+p_8\xi_2)^m= (p_3
\xi_0+p_4\xi_1+p_5\xi_2)^n$$ 
(taking $(1:x:y)$ for coordinates in $\P^2$). Pulling back through
$\loc\pi$ (in the first case written above; the second case is
analogous) gives
$$(\xi_0+s_1\xi_1+s_2\xi_2)^{n-m}(\xi_0+s_7\xi_1+s_8\xi_2)^m s_3^{mn}
s_6^{Bm}= (\xi_0+s_4\xi_1+s_5\xi_2)^n s_3^{mn} s_6^{An}$$
Using $Bm-An=1$, we see that the proper transform of this
point@-condition in $\loc{\Til V}$ has equation
$$(\xi_0+s_1\xi_1+s_2\xi_2)^{n-m}(\xi_0+s_7\xi_1+s_8\xi_2)^m s_6=
(\xi_0+s_4\xi_1+s_5\xi_2)^n $$
Note that $s_3$ does not appear in this equation. Next, recall that
the support of the base locus of $c$ in $\P^8$ is $\P^2\times C$; as
$C$ is parametrized by $(t^m,t^n)$, we may parametrize the (affine
part of the) support of the base locus by
$$(k_1,k_2,t)
\mapsto \pmatrix 1 & k_1 & k_2 \\
t^m & t^m k_1 & t^m k_2\\
t^n & t^n k_1 & t^n k_2
\endpmatrix$$
Using $Bm-An=1$ again we lift this to
$$(k_1,k_2,t) \mapsto \pmatrix
\boxed{\loc{\Til V}} & k_1 & k_2\\
t & k_1 & k_2\\
1 & k_1 & k_2
\endpmatrix$$
and observe that a point of this subvariety of $\loc{\Til V}$ lies
above the special point of $C$ $\iff t=0\iff s_3=0$. Since the
equations of the point@-conditions do not involve $s_3$, their behavior
over such a point is the same as over any point with nonzero $t$. As
no point with $t\ne 0$ is a flex on $C$, we know from \cite{A-F1},
Proposition~2.7, that two `global' blow@-ups resolve the
indeterminacies of $c$ at such points.\qed\enddemo

As stated in the proof of the theorem, the two blow@-ups needed to
remove the indeterminacies of $\loc c$ are the two blow@-ups,
discussed in \cite{A-F1}, resolving the map over nonsingular non@-flex
points of $C$. We refer the reader to \cite{A-F1} for a more thorough
description of the centers of these blow@-ups, and freely use that
information in \S3. Here we will just recall that the first `global'
blow@-up will have a nonsingular irreducible 3@-dimensional center
(the proper transform of $\P^2\times C$ in $\loc{\Til V}$); after
blowing up this locus, the point@-conditions meet along a
4@-dimensional locus, in fact a $\P^1$@-bundle over the preceding
center. The point@-conditions are separated from each other by blowing
up this last locus. This $\P^1$ bundle is described in the discussion
preceding Proposition~2.2 in \cite{A-F1}.

Now we move to the coordinate description of $\loc{\Til V}$ used
above. All the varieties we consider are obtained by a sequence of
blow@-ups over $\P^8$, and 
inherit a right action of 
$\PGL(3)$ since the centers of the blow@-ups are invariant.
We say that a chart in any such variety is
{\it essential\/} if every point of the variety can be moved to that
chart by this action.

\proclaim{Proposition 2.6} The variety $\loc{\Til V}$ admits an
essential chart with the coordinate description specified
above.\endproclaim

\demo{Proof} Let $\loc V_1$ be the $\ell_1$@-directed
blow@-up of $\P^8$ along $B$ in the direction of $P$. In order to
study this variety, we first obtain local coordinates for the base
locus of $c$.
Writing out the matrix with kernel on $x_0+k_1 x_1+k_2 x_2=0$ and image
$(1:t^m:t^n)$ gives a local parametrization
$$(k_1,k_2,t) \mapsto \pmatrix 1 & k_1 & k_2 \\
t^m & t^m k_1 & t^m k_2\\
t^n & t^n k_1 & t^n k_2
\endpmatrix$$
for $\P^2\times C$. Setting $t=0$ selects the distinguished $\P^2$,
locally parametrized by
$$(k_1,k_2,0) \mapsto \pmatrix 1 & k_1 & k_2 \\
0 & 0 & 0\\
0 & 0 & 0
\endpmatrix$$
So $B$ has equations
$$p_3=p_4=p_5=p_6=p_7=p_8=0$$
in our chart. As for $P$, the matrices with image contained in the
line $y=0$ are in the form
$$\pmatrix 
* & * & *\\
* & * & *\\
0 & 0 & 0
\endpmatrix$$
so $P$ has equations
$$p_6=p_7=p_8=0$$
The distinguished 4@-fold in $\loc V_1$ is the intersection $\loc B_1$
of the proper transform of $P$ and the last exceptional divisor $\loc
E_1$ of the sequence producing the directed blow@-up. In fact we
can also consider the $j$@-directed blow@-up for all $1\le j<\ell_1$,
and a simple inductive computation shows that at each stage the proper
transforms of the point@-conditions meet along the proper transform of
$\P^2\times C$, and along the locus $B_i$ obtained by intersecting the
proper transform of $P$ with the last exceptional divisor. If $e>1$,
the same holds for the $\loc V_1$ (as we will see below), so we only
need to examine $\loc V_1$ near $\loc B_1$. Now, by part (4) of
Lemma~2.2,
a neighborhood of $\loc B_1$ is covered by charts $U_3$, $U_4$,
$U_5$ with local parameters $(q_1,\dots,q_8)$ so that the map $U_j
@>>> \P^8$ is given by
$$p_i=\left\{\aligned
q_i  \quad\quad &i=j, i=1,2\\
q_i q_j \quad\quad &i=3,4,5\,\text{ but $i\ne j$}\\
q_i q_j^{\ell_1}\quad\quad &i=6,7,8\\
\endaligned\right.$$
\proclaim{Claim} The chart $U_3$ is essential.\endproclaim
That is, we claim that we can use the action of $\PGL(3)$ to move
points from the other charts to this chart. The proof of this fact is
a simple but tedious coordinate computation, which we leave to the
reader.

The consequence of the claim is that it is not restrictive to choose
local coordinates $q_i$ on $\loc V_1$ so that the blow@-up map $\loc V_1
@>>> \P^8$ is given by
$$\pmatrix
1 & p_1 & p_2\\
p_3 & p_4 & p_5\\
p_6 & p_7 & p_8
\endpmatrix=\pmatrix
\boxed{\loc V_1} & q_1 & q_2\\
q_3 & q_3 q_4 & q_3 q_5\\
q_3^{\ell_1} q_6 & q_3^{\ell_1} q_7 & q_3^{\ell_1} q_8
\endpmatrix$$
The equation of the exceptional divisor $\loc E_1$ is $q_3=0$ in these
coordinates. The equation of the point@-condition corresponding to
$(\xi_0:\xi_1:\xi_2)$ in $\P^8$ is
$$(\xi_0+p_1\xi_1+p_2\xi_2)^{n-m}(p_6\xi_0+p_7\xi_1+p_8\xi_2)^m= (p_3
\xi_0+p_4\xi_1+p_5\xi_2)^n\quad,$$ 
which pulled back via the above map gives
$$(\xi_0+q_1\xi_1+q_2\xi_2)^{n-m}(q_6\xi_0+q_7\xi_1+q_8\xi_2)^m
q_3^{\ell_1 m}= (\xi_0+q_4\xi_1+q_5\xi_2)^nq_3^n\quad;$$
clearing a common factor of $q_3^{\ell_1 m}$ (notice $n-\ell_1
m=m_2\ge 0$), we obtain the equation of the point@-condition in $\loc
V_1$:
$$(\xi_0+q_1\xi_1+q_2\xi_2)^{n-m}(q_6\xi_0+q_7\xi_1+q_8\xi_2)^m=
(\xi_0+q_4\xi_1+q_5\xi_2)^nq_3^{m_2}\quad.$$
If $m_2=0$, that is $e=1$, then $\loc{\Til V}=\loc V_1$; setting
$s_i=q_i$, $A=0$, and $B=1$ gives the prescribed coordinate
description, and we are done in this case.

Otherwise, this shows that along $\loc E_1$ (i.e., setting $q_3=0$)
the point@-conditions meet along the locus with equations
$$q_6=q_7=q_8=0$$
that is, the intersection of the proper transform of $P$ with $\loc
E_1$. This is the distinguished 4@-fold, $\loc B_1$.\vskip 6pt

Next, we perform the $\ell_2$@-directed blow@-up of $\loc V_1$ along
$\loc B_1$, in the direction of $\loc E_1$, obtaining $\loc V_2$, with
exceptional divisor $\loc E_2$. The discussion is very similar to the
discussion of the first step; now $\loc B_1\subset \loc E_1$ are given
by the ideals $(q_3,q_6,q_7,q_8)\supset (q_3)$, so by Lemma~2.2
we can choose a chart in $\loc V_2$ with coordinates $r_i$, such that
the map $\loc V_2 @>>> \loc V_1$ is given by
$$\pmatrix
\boxed{\loc V_1} & q_1 & q_2\\
q_3 & q_4 & q_5\\
q_6 & q_7 & q_8
\endpmatrix=\pmatrix
\boxed{\loc V_2} & r_1 & r_2\\
r_3 r_6^{\ell_2}& r_4 & r_5\\
r_6 & r_6 r_7 & r_6 r_8
\endpmatrix\quad.$$
Again, the reader should have no difficulties checking that this chart
is essential. 

The new exceptional divisor $\loc E_2$ is given by $r_6=0$ in these
coordinates. Pulling back the point@-conditions from $\loc V_1$ and
clearing a common factor of $r_6^{\ell_2 m_2}$ shows that the equation
of the point@-condition corresponding to $(\xi_0:\xi_1:\xi_2)$ in
$\loc{V_2}$ is
$$(\xi_0+r_1\xi_1+r_2\xi_2)^{n-m}(\xi_0+r_7\xi_1+r_8\xi_2)^mr_6^{m_3}=
(\xi_0+r_4\xi_1+r_5\xi_2)^nr_3^{m_2}\quad.$$
If $m_3=0$, that is $e=2$, then $\loc{\Til V}=\loc V_2$, and we have
reached the desired coordinate expression.

Otherwise, we see that the intersection of all point@-conditions along
$\loc E_2$ is the locus with equations
$$r_3=r_6=0\quad,$$
giving the distinguished 6@-fold $\loc B_2$.\vskip 6pt

Having reached this stage, the expression for the point@-condition
$$(\xi_0+r_1\xi_1+r_2\xi_2)^{n-m}(\xi_0+r_7\xi_1+r_8\xi_2)^mr_6^{m_3}=
(\xi_0+r_4\xi_1+r_5\xi_2)^nr_3^{m_2}$$
is so symmetric that the remaining blow@-ups can be all understood
together. 
Assuming that we have defined $\loc V_i$, we will have

---either $i$ odd, equation of $\loc E_i$: $r_3=0$; and equation
of the point@-conditions
$$(\xi_0+r_1\xi_1+r_2\xi_2)^{n-m}(\xi_0+r_7\xi_1+r_8\xi_2)^m
r_6^{m_i}=(\xi_0+r_4\xi_1+r_5\xi_2)^nr_3^{m_{i+1}}\quad;$$

---or $i$ even, equation of $\loc E_i$: $r_6=0$; and equation
of the point@-conditions
$$(\xi_0+r_1\xi_1+r_2\xi_2)^{n-m}(\xi_0+r_7\xi_1+r_8\xi_2)^m
r_6^{m_{i+1}}=(\xi_0+r_4\xi_1+r_5\xi_2)^nr_3^{m_i}\quad.$$

As long as $m_{i+1}>0$, the point@-conditions meet on $\loc E_i$ along
the 6@-fold $\loc B_i$ defined by $r_3=r_6=0$ (in both
cases). Applying again part (4) of Lemma~2.2, we see that the
$\ell_{i+1}$@-directed blow@-up of $\loc V_i$ along $\loc B_i$ in the
direction of $\loc E_i$ produces a $\loc V_{i+1}$ with the data
prescribed above;
in particular, we see that this automatically chooses the essential
chart in each successive blow@-up.
Notice in passing that at each stage $\loc B_i$ is
the intersection of $\loc E_i$ with the proper transform $\loc {\Til
E}_{i-1}$; and the two divisors are swapped from one stage to the
next. In particular, the restriction of $\loc {\Til E}_{i-1}$ to $\loc
B_i\cong \loc B_{i+1}$ equals the restriction of $\loc E_{i+1}$. This
fact will be used in \S3.

At the $e$@-th stage we will have $m_e=gcd(m,n)=1$ and $m_{e+1}=0$, so
the point@-conditions will have equation
$$(\xi_0+r_1\xi_1+r_2\xi_2)^{n-m}(\xi_0+r_7\xi_1+r_8\xi_2)^mr_6
=(\xi_0+r_4\xi_1+r_5\xi_2)^n$$
for odd $e$, and
$$(\xi_0+r_1\xi_1+r_2\xi_2)^{n-m}(\xi_0+r_7\xi_1+r_8\xi_2)^m
=(\xi_0+r_4\xi_1+r_5\xi_2)^nr_3$$
for even $e$. Writing the map to $\P^8$ explicitly shows that this
gives the claimed coordinate description of $\loc V_e=\loc {\Til V}$,
as needed.\qed\enddemo

It follows from the explicit equations obtained in this proof that
the multiplicity of the point@-conditions along the various centers of
blow@-up also mirrors the multiplicity of $C$ at the centers of the
blow@-ups resolving it. So this multiplicity is $m_1$ for the $\ell_1$
blow@-ups giving the first directed blow@-up, $m_2$ for the second
batch, etc. This information will be used in \S3.

We will also need the multiplicities of the $\P^7$'s obtained as
point@-conditions for the lines of the basic triangle, so we note here
that these also mirror the corresponding multiplicities of the lines
in the blow@-ups resolving the curve. Explicitly, for the blow@-ups
examined here 

---a $\P^7$ corresponding to the line connecting the two cusps of
$C$ ($\lambda$ in the notation of \S1) has multiplicity~1 along the
first center $B$ of the first directed blow@-up, and multiplicity~0
at all other centers;

---a $\P^7$ corresponding to the line $\mu$ supporting the tangent
cone to $C$ at the cusp under consideration has multiplicity~1 along
all the centers of the blow@-ups giving the first directed blow@-up;
multiplicity~1 along the first center $\loc B_1$ of the second
directed blow@-up; and multiplicity~0 at all other centers;

---a $\P^7$ corresponding to the line $\overline \mu$ at infinity has
multiplicity~0 along all centers.

\subheading{\S2.3. Blow@-ups for the general case}
It is now a simple matter to go from the case of one curve, treated in
\S2.2, to the case of many. Again, the more general curves $C$ we
consider in this paper are arbitrary unions (with multiplicities) of
elements of the pencil
$$x^n=\alpha y^m z^{n-m}\quad (\alpha\ne 0)$$
together with multiples of the lines $\lambda, \mu, \overline\mu$ of
the basic triangle. As pointed out in \S2.2, removing the
indeterminacies of the corresponding rational map $c$ amounts to
separating the point@-conditions; so we have to understand what the
point@-conditions of $C$ look like, and how they behave under the
blow@-ups described in \S2.2.
\proclaim{Proposition 2.7} For a curve $C$ as above, the
point@-condition in $\P^8$ corresponding to a point $p\in \P^2$
consists of the union of the point@-conditions of each component,
each appearing with multiplicity equal to the multiplicity of the
corresponding component.\endproclaim
This should be clear: if $F=0$ is an equation for $C$, then the
equation of the point@-condition corresponding to $p$ is the vanishing
of 
$$F(\varphi(p))=0\quad.$$
This polynomial (in $\varphi$) factors according to how $F$ factors.

The supports of the point@-conditions of $C$ are therefore unions of
point@-conditions considered in \S2.2 (for different $\alpha$'s),
and of copies of the three $\P^7$'s corresponding to $\lambda$, $\mu$,
$\overline\mu$ mentioned at the end of \S2.2. 

Disregarding the lines of the basic triangles for a moment, note that
different irreducible curves from the same pencil as above have the
same history through this blow@-up sequence. The situation in the
plane mirrors precisely the situation at the level of
point@-conditions: different curves determine the same centers $B$,
$\loc B_i$, and the corresponding point@-conditions have the same
multiplicities along these loci. Further, the curves are separated at
the very last stage, and correspondingly the base locus of the lifted
map $\loc {\Til V} \dashrightarrow \P^N$ consists of the disjoint
union of copies of $\Til C\times\P^2$ (where $\Til C$ denotes the
normalization of a single curve of type $(m,n)$). For each of these,
Theorem~2.4 shows that two (`global') blow@-ups will suffice
to remove the indeterminacies.

In other words, the same sequence of local blow@-ups used for one
curve of type $(m,n)$, followed by two global blow@-ups for each
component, removes the indeterminacies for any finite union of such
curves. The multiplicities with which these appear are irrelevant to
this discussion.

To account for the lines in the basic triangle, we need to keep
track of the three pencils of hyperplanes of $\P^8$ corresponding to
the points on these three lines. The relevant data is
implicit in the multiplicity statement at the end of \S2.2:

---the $\P^7$'s corresponding to $\lambda$ are separated from the
point@-conditions corresponding to curves of type $(m,n)$ after the
first blow@-up of the sequence giving the first directed blow@-up;

---the $\P^7$'s corresponding to $\mu$ are separated from the
point@-conditions corresponding to curves of type $(m,n)$ at the first
blow@-up of the sequence giving the second directed blow@-up over the cusp
$\lambda\cap\mu$;

---the $\P^7$'s corresponding to $\overline\mu$ are separated from the
point@-conditions corresponding to curves of type $(m,n)$ at the first
blow@-up of the sequence giving the second directed blow@-up over the
cusp $\lambda\cap\overline\mu$.

That is, after the local blow@-ups of \S2.2 have been performed over
both cusps, and after the two global blow@-ups of \S2.2 have removed
indeterminacies arising from the type@-$(m,n)$ components of $C$, we
still have three groups of hypersurfaces, corresponding to the three
lines of the triangle. Again, it is easily checked that the incidence
of these groups of hypersurfaces reflects the incidence of the
corresponding proper transforms of the lines in the plane:
$$\epsffile{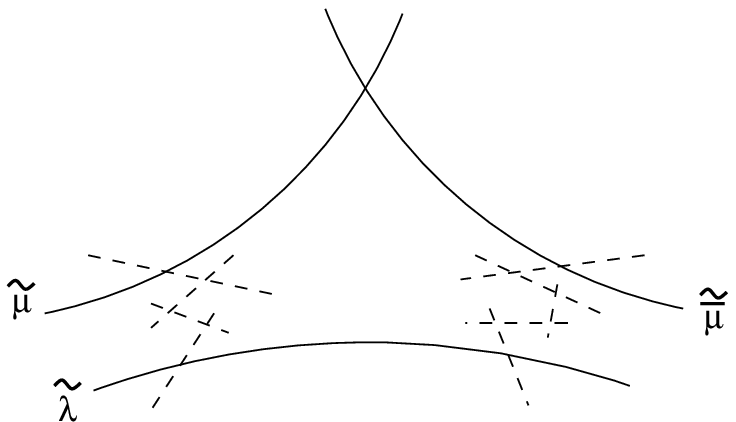}$$

Finally we deal with these hypersurfaces. The intersection of the
hypersurfaces in each group is a five@-dimensional variety (the proper
transform of the $\P^5$ of rank@-2 matrices with image the
corresponding line of the triangle); further, the two
five@-dimensional varieties corresponding to $\mu$ and $\overline \mu$
still meet along a $\P^2$, corresponding to rank@-1 matrices whose
image is the point of intersection $\mu\cap\overline\mu$.

By our good luck (as the reader can see by performing the relevant
computation, using the coordinates given in \S2.2), the `obvious'
strategy works: blowing up along the $\P^2$ corresponding to
$\mu\cap\overline\mu$, and then along the proper transforms of the
three $\P^5$ described above, finally produces a variety $\Til V$
satisfying the condition in the statement of Theorem~2.1.

\head \S3. Degree computations\endhead

With the coordinate analysis of \S2~behind us, we are ready to set
up the intersection theoretic part of the computation. The discussion
leading to Proposition~3.2 below reduces the computation of the
degree of an orbit closure to the computation of an intersection
product of divisors on the variety $\Til V$ we constructed in \S2. Our
main tool will then be Proposition~2.3, by which we keep track of
the intersection of divisors under directed blow@-ups.

The information needed to apply this formula consists of the
multiplicity of the divisors at the center of blow@-up, together with
the Chern classes of the relevant normal bundles. The first piece of
information is listed at the end of \S2.2; the second will be obtained
along the way, mostly by using part (3) of Lemma~2.2.

Here is the main reduction.
We want to compute the
degree of the orbit closure $\overline{\Cal O}_C\subset \P^N$,
assuming this has dimension 7 (which is the case for most choices of
the parameters $m,n,$ etc.). In the set@-up of \S2, we have obtained a
completion $\Til V$ of $\PGL(3)$ over which the action extends to a
regular map to $\P^N$; $\Til V$ was obtained by suitably blowing up
$\P^8$:
$$\diagram
{\PGL(3)\subset \Til V\phantom{\PGL(3)\subset}} \dto_\pi \drto^{\tilde c}\\
{\P^8} \xdashed[0,1]^{c}|>\tip & {\P^N} \\
\enddiagram$$
This realizes $\overline{\Cal O}_C$ as the image of $\Til V$ by
$\tilde c$.

Now, seven general hyperplanes intersect $\overline{\Cal O}_C$
transversally at $\deg\overline{\Cal O}_C$ points of $\Cal O_C$. The
inverse image of these points in $\PGL(3)$ will be $\deg\overline{\Cal
O}_C$ translated copies of the stabilizer of $C$. It follows that
$$({\tilde c}^*h)^7=\deg\overline{\Cal O}_C\,[Z]\tag *$$
where $h$ is the class of a hyperplane in $\P^N$, and $Z$ is the cycle
obtained by closing up in $\Til V$ the stabilizer of $C$.

By construction, the class ${\tilde c}^*h$ is represented by a
`point@-condition' $\Til W\subset\Til V$; that is, by the proper
transform of the hypersurface $W\subset \P^8$ consisting of matrices
$\phi$ mapping a fixed point $p\in \P^2$ to a point of $C$.
Note that $\pi(Z)$ consists of the closure in $\P^8$ of the stabilizer
of $C$. As mentioned after the statement of Theorem~1.1,
the number of components of the stabilizer of $C$ depends on
symmetries of the specific $S$@-tuple of points in $\A^1$
corresponding to the non@-linear components of $C$. Explicitly, assume
that $C$ is given by the equation
$$
x^r y^q z^{\overline q} \prod_i (x^n-\alpha_i y^m z^{\overline m})
^{s_i}=0\quad,
$$
so that it corresponds to the $S$@-tuple in $\A^1_{\alpha}$ given by
the equation
$$\prod_i (\alpha-\alpha_i)^{s_i}=0\quad.$$
The components of the stabilizer of $C$ depend on automorphisms $\A^1
@>>> \A^1$ fixing this $S$@-tuple. The precise statement (whose proof
is left to the reader) is

\proclaim{Lemma 3.1} With notations as above, assume that the orbit of
$C$ has dimension 7.
Then if $n\ne 2$ or $q\ne \overline q$ the number $A$ of components of
the stabilizer of $C$ equals the number of automorphisms $\A^1 @>>> \A^1$,
$\alpha\mapsto u\alpha$ (with $u$ a root of unity) preserving the
$S$@-tuple corresponding to $C$; when $n=2$ and $q={\overline q}$, $A$
equals twice this number.
\endproclaim

The extra automorphisms in the latter case come from the switch
$y\leftrightarrow z$. 

\example{Examples} (1) $C$ given by $(x^3-yz^2)(x^3-2yz^2)=0$. The
corresponding $S$@-tuple is given by $(\alpha-1)(\alpha-2)=0$, that is,
by the pair of points $\alpha=1$, $\alpha=2$ in $\A^1$. The only
automorphism of $\A^1$ of the kind specified in the statement and
preserving this pair of points is the identity, so $A=1$ in this case
(as in most others).

(2) $C$ given by $(x^2-yz)(x^2+yz)=0$. The corresponding $S$@-tuple is
$\alpha=\pm 1$. Two automorphisms preserve this pair: the identity and
$\alpha\mapsto -\alpha$. Since $n=2$ and $q=\overline q=1$, we have
$A=4$ in this case.
\endexample

It is easily checked that each component of the closure of the
stabilizer is a copy of the curve $x^n=y^m z^{\overline m}$. For
example, the identity component of the stabilizer of $C$ consists of
the diagonal matrices with entries $(1,t^{\overline m},t^{-m})$; since
$(m,\overline m)=1$, its closure has equation  $1=p_4^m p_8^{\overline
m}$ in the coordinates of \S2. In particular, the degree of $\pi(Z)$
equals $An$.

Pushing forward (*) to $\P^8$ and intersecting by a general hyperplane
$H$, we see then that
$$\int H\cdot \pi_*({\Til W}^7)=A\,n\,\deg\overline{\Cal O}_C\quad.$$

Observing that the inverse image $\pi^{-1} H$ of a general hyperplane
equals its proper transform $\Til H$, and applying the projection
formula, we conclude:
\proclaim{Proposition 3.2} If $\dim\overline{\Cal O}_C=7$, then
$$\deg\overline{\Cal O}_C=\frac 1{A\,n}\int_{\P^8} \pi_*(\Til H\cdot \Til
W^7)\quad.$$
\endproclaim

Our goal is therefore to perform the intersection product on the
right@-hand@-side of this formula.

The reader now sees why we stated the formula in Proposition~2.3,
comparing intersection products of divisors and of their
proper transforms under directed blow@-ups. The role of the different
divisors considered in that formula might not be immediately apparent,
however, and the next lemma should clarify it. We denote by $H$ the
general hyperplane in $\P^8$; by $H_\lambda$, $H_\mu$, $H_{\overline\mu}$
respectively hyperplanes obtained as point@-conditions relative to
the lines $\lambda$, $\mu$, $\overline\mu$. Further, we denote by $X$
a point@-condition in $\P^8$ relative to the part of $C$ consisting of
the union of type@-$(m,n)$ curves.
\proclaim{Lemma 3.3} With these notations,
$$A\,n\,\deg\overline{\Cal O}_C=7!\sum
\frac{q^{j_\mu}}{j_\mu!}\, \frac{{\overline q}^{j_{\overline\mu}}}
{j_{\overline\mu}!}\, \frac{r^{j_\lambda}} {j_\lambda!}\,\frac{1}
{j_c!} \int \pi_*\left(\Til H\cdot \Til H_\mu^{j_\mu}\cdot \Til
H_{\overline\mu}^{j_{\overline\mu}} \cdot \Til H_\lambda^{j_\lambda}
\cdot \Til X^{j_c}\right)\quad,$$
where the summation runs over all $0\le j_\mu\le 2$, $0\le
j_{\overline\mu}\le 2$, $0\le j_\lambda\le 2$, $0\le j_c\le 7$ such
that $j_\mu+j_{\overline\mu}+j_\lambda+j_c=7$.\endproclaim
\demo{Proof} As observed in Proposition~2.7,
point@-conditions of a reducible curve split into the point@-conditions
of its irreducible components. This implies
$$\Til W=q\Til H_\mu+\overline q\Til H_{\overline\mu}+r \Til
H_\lambda+\Til X\quad.$$
The formula follows then from Proposition~3.2, once one
observes that if $\ell$ is any of $\lambda$, $\mu$, $\overline\mu$,
then $\Til H_\ell^3=0$: indeed, a line does not contain three general
points, so the intersection of three point@-conditions of a line must
be empty.\qed\enddemo

By this lemma, we are reduced to computing intersection products
$$\Til H\cdot \Til H_\mu^{j_\mu}\cdot \Til
H_{\overline\mu}^{j_{\overline\mu}} \cdot \Til H_\lambda^{j_\lambda}
\cdot \Til X^{j_c}\quad,$$
for $0\le j_\mu\le 2$, $0\le j_{\overline\mu}\le 2$, $0\le
j_\lambda\le 2$, $0\le j_c\le 7$ such that
$j_\mu+j_{\overline\mu}+j_\lambda+j_c=7$. The divisors $H_\lambda$,
$H_\mu$, $H_{\overline\mu}$ will play the role of the divisors `of
type $Y_i$' in Proposition~2.3.

\subheading{\S3.1. Local blow@-ups} Now we move to the core of the
computation. Proposition~2.3 will be used iteratively to evaluate
the intersection product listed above on successively higher and
higher level blow@-ups. At each directed blow@-up, the formula
evaluates a correction term measuring by how much the intersection
product changes upon taking proper transforms. The starting point is
the intersection product in $\P^8$,
$$H\cdot H_\mu^{j_\mu}\cdot H_{\overline\mu}^{j_{\overline\mu}} \cdot
H_\lambda^{j_\lambda} \cdot X^{j_c}\quad:$$
since (with the notation of \S1) $X$ has degree $Sn$, this is simply
$$(Sn)^{j_c}$$
by B\'ezout's Theorem. Summing up as in Lemma~3.3, we get:
$$7!\left(\sum_{0\le j_\mu\le 2,\, 0\le j_{\overline\mu}\le 2,\,  0\le
j_\lambda\le 2,\, j_c=7-j_\mu-j_{\overline\mu}-j_\lambda}
\frac{q^{j_\mu}}{j_\mu!}\, \frac{{\overline q}^{j_{\overline\mu}}}
{j_{\overline\mu}!}\, \frac{r^{j_\lambda}} {j_\lambda!}\,\frac{(S
n)^{j_c}} {j_c!}\right)\quad.$$
This unpleasant expression prompts us to establish the following:
\remark{Convention} We are going to treat the multiplicities $q$,
$\overline q$, $r$, $S$ as variables, and impose that $q^3=\overline
q^3=r^3=0$.\endremark
This takes care automatically of the bounds for the $j$'s in the
summation, so that the B\'ezout term simply becomes
$$(Sn+r+q+\overline q)^7\quad.$$
The geometric reason behind the convention is that the
self@-intersection of three or more point@-conditions in $\Til V$
corresponding to lines must vanish, as was mentioned above. Imposing
this from the start saves us some computational time: in practice, all
the terms that we discard at this stage would be cancelled anyway along the
blow@-up process, so we can ignore them. The important caveat to keep
in mind is that one {\it may not substitute the multiplicities for
their value before expanding expressions in which they appear.\/} All
such expressions must be expanded, and the relations $q^3=\overline
q^3=r^3=0$ must be applied, before any substitution can be made.\vskip
6pt

Next, we deal with the correction term due to the $e$ directed
blow@-ups over the cusp at $\lambda\cap\mu$~(where $e=$ number of
lines in the Euclidean algorithm for $(m,n)$). By symmetry, we will
get a similar contribution for the cusp $\lambda\cap\overline\mu$. In
the next subsection we will evaluate analogous contributions due to
the other (`global') blow@-ups.

Recall our notation: the Euclidean algorithm performed on $m=m_1$ and $n$
gives
$$\align
n &=\ell_1 m_1+m_2\\
m_1 &=\ell_2 m_2+m_3\\
&\cdots\\
m_{e-1}&=\ell_e m_e
\endalign$$
with all $m_i$, $\ell_i$ positive
integers, $0<m_{i+1}<m_i$, and $m_e=gcd(m,n)$; we are in fact assuming
$m_e=1$. The first blow@-up from Proposition 2.6~is
the $\ell_1$@-directed blow@-up of $\P^8$ along $\P^2$ in the
direction of $\P^5$, where these subspaces are defined immediately
preceding the statement of Proposition~2.6. The multiplicities
of the (proper transforms of the) divisors we need to intersect were
discussed at the end of \S2.2,
and are as follows:

---for the general hyperplane: 0 for all centers of the $\ell_1$ blow@-ups;

---for the $j_\mu$ hyperplanes corresponding to the tangent cone
to $C$ at the point: 1 for all centers;

---for the $j_{\overline\mu}$ hyperplanes corresponding to the tangent
cone to $C$ at the other point: 0 for all centers;

---for the $j_\lambda$ hyperplanes corresponding to the line joining
the two distinguished points of $C$: 1 for the first center, 0 for the
remaining $\ell_1-1$;

---for the $j_c$ point@-conditions $X$: $S m_1$ at all centers.

\noindent (indeed, each support of a component of $X$ has multiplicity
$m_1$, and $S=$the sum of the multiplicities of the components).
Further, $X$ has degree $S n$. Also observe that we have $j_\lambda \le
2<\codim_{\P^2} \P^5$ terms with `mixed multiplicities', as is
necessary in order to apply Proposition~2.3.

Finally, denoting by $k$ the hyperplane class in
$\P^2$, we have
$$c(N_{\P^2}\P^5)=(1+k)^3\quad,\quad c(N_{\P^5}\P^8)=(1+k)^3$$
and therefore
$$c(N_{\P^2}^{(\ell_1)}\P^5)=(1+\ell_1 k)^3$$
and Proposition~2.3 evaluates the correction term due to the
first directed blow@-up:
$$\ell_1^{3-j_\lambda} \int_{\P^2} \frac {k (1+k)^{j_\mu} (1+\ell_1
k)^{j_\lambda} k^{j_{\overline \mu}} (S m_1+S n k)^{j_c}}{(1+\ell_1
k)^3 (1+k)^3}\cap[\P^2]$$
that is, with minimal manipulations:
$$\boxed{
S^{j_c}\ell_1^{3-j_\lambda}\int_{\P^2} \frac {k^{j_{\overline\mu}+1}
(m_1+n k)^{j_c}} {(1+\ell_1 k)^{3-j_\lambda} (1+k)^{3-j_\mu}}\cap[\P^2]}$$

We will see that, remarkably, {\it this box contains all the
information necessary to compute the `local' contributions.\/}
However, to understand this we have to write similar terms for the
other directed blow@-ups.\vskip 6pt

For the $\ell_2$@-directed blow@-up of $\loc V_1$ along $\loc B_1$
in the direction of $\loc E_1$, the formula will evaluate the term as
a degree over the 4@-fold $\loc B_1$. As $\loc B_1$ was obtained as
the distinguished $\Til B$ in the $\ell_1$@-directed blow@-up
considered at the first stage, points (1) and (3) in Lemma~2.2 
give that $\loc B_1$ is the projectivization of the normal
bundle of $\P^2$ in $\P^5$, and that
$$c(N_{\loc B_1}\loc E_1)=c(N_{\P^5}\P^8\otimes \Cal O(-\ell_1)) =
(1+k-\ell_1 e_1)^3$$
where $k, e_1$ are respectively the pull@-back of $k$ from $\P^2$ and
the restriction of the class of $\loc E_1$ to $\loc B_1$. The
description of $\loc B_1$ as $\P(N_{\P^2}\P^5)$ gives easily all the
information needed to perform computations in the intersection ring of
$\loc B_1$: to evaluate explicitly the term we are going to write in a
moment we would push@-forward to $\P^2$, then use
$$e_1^4\cdot [\loc B_1] \mapsto 6 k^2\cdot [\P^2], \quad e_1^3
\cdot [\loc B_1] \mapsto 3 k\cdot [\P^2],\quad e_1^2\cdot [\loc
B_1] \mapsto [\P^2],\quad e_1\cdot [\loc B_1] \mapsto 0$$
(this follows immediately from $\frac {[\P(N_{\P^2}\P^5)]}
{(1+e_1)} \mapsto c(N_{\P^2}{\P^5})^{-1}\cap [\P^2]=\frac{[\P^2]}
{(1+k)^3}$, cf.~\cite{Fulton}).

Next, we list here the multiplicities of the divisors and the classes of
their pull@-backs to $\loc B_1$:

---for the general hyperplane: multiplicity 0, class $k$;

---for the $j_\lambda$ hyperplanes corresponding to the line joining
the two distinguished points of $C$: 0 for all centers; class
$=k-e_1$;

---for the $j_\mu$ hyperplanes corresponding to the tangent cone:
1 along the first center, 0 along the remaining $\ell_2-1$; class
$=k-\ell_1 e_1$;

---for the $j_{\overline\mu}$ hyperplanes corresponding to the tangent
cone to $C$ at the other point: 0 for all centers; class $=k$;

---for the $j_c$ point@-conditions: $S m_2$ along all centers; class
$=S n k-S \ell_1 m_1 e_1$.

\noindent Note that this time the number of divisors with mixed
multiplicity is $j_\mu\le 2<\codim_{\loc B_1}\loc E_1$, as needed to
apply Proposition~2.3.

Finally, from the above:
$$c({N_{\loc B_1}^{(\ell_2)}{\loc E_1}})=(1+\ell_2 k-\ell_1\ell_2 e_1)^3
\quad,\quad c(N_{\loc E_1}\loc V_1)=(1+e_1)$$
and we are ready to apply Proposition~2.3, which gives
$$\boxed{\ell_2^{3-j_\mu} \int_{\loc B_1} \frac {k 
(1+\ell_2 k-\ell_1 \ell_2 e_1)^{j_\mu} (k-e_1)^{j_\lambda}
k^{j_{\overline\mu}} (S m_2+S n k-S \ell_1 m_1 e_1)^{j_c}} {(1+\ell_2
k-\ell_1\ell_2 e_1)^3(1+e_1)}\cap[\loc B_1]}$$

Again, this could be somewhat simplified.

The remaining blow@-ups all have isomorphic centers, so we can
describe the corresponding terms uniformly. Recall from
Proposition~2.6 that
at the $i$@-th stage ($i>2$) we are performing the
$\ell_i$@-directed blow@-up of $\loc V_{i-1}$ along $\loc B_{i-1}$ in
the direction of $\loc E_{i-1}$; here $\loc B_{i-1}$ is a 6@-fold, and
again we can describe it concretely by using Lemma~2.2, point (1):
it is the projectivization of the normal bundle to $\loc B_1$ in $\loc
E_1$. An alternative description is as the intersection of the proper
transform of $\loc E_{i-2}$ with $\loc E_{i-1}$; denoting by $e_j$ the
restriction of $\loc E_j$, the proper transform of $\loc E_{i-2}$ will
restrict to $e_{i-2}-\ell_{i-1} e_{i-1}$, so
$$c(N_{\loc B_{i-1}}\loc E_{i-1})=(1+e_{i-2}-\ell_{i-1} e_{i-1})$$
Also notice that by the interchange of exceptional divisors we
observed toward the end of the proof of Proposition~2.6
we have
$$e_i=e_{i-2}-\ell_{i-1} e_{i-1}$$

Multiplicities and class of the divisors:

---for the general hyperplane: multiplicity 0, class $k$;

---for the $j_\lambda$ hyperplanes corresponding to the line joining
the two distinguished points of $C$: multiplicity 0, class $=k-e_1$;

---for the $j_\mu$ hyperplanes corresponding to the tangent cone:
multiplicity 0, class $=k-\ell_1 e_1-e_2$;

---for the $j_{\overline\mu}$ hyperplanes corresponding to the tangent
cone to $C$ at the other point: multiplicity 0; class $=k$;

---for the $j_c$ point@-conditions $\Til X$: $S m_i$ along all
centers; class $=S n k-S \ell_1 m_1 e_1-\dots -S \ell_{i-1}
m_{i-1} e_{i-1}$.

{From} the above, the classes of the relevant bundles are
$$c({N_{\loc B_{i-1}}^{(\ell_i)}{\loc E_{i-1}}})=(1+\ell_ie_{i-2}-
\ell_{i-1} \ell_i e_{i-1})\quad,\quad c(N_{\loc E_{i-1}}\loc V_{i-1})
=(1+e_{i-1})$$
and Proposition~2.3 evaluates the $i$@-th term:
$$\ell_i \int_{\loc B_{i-1}} \frac {k (k-\ell_1 e_1-e_2)^{j_\mu}
(k-e_1)^{j_\lambda} k^{j_{\overline\mu}} (S\cdot\text{p.c.})^{j_c}}
{(1+\ell_i e_{i-2}-\ell_{i-1}\ell_i e_{i-1})(1+e_{i-1})}\cap[\loc B_{i-1}]$$
where
$$\text{p.c.}=m_i+n k- \ell_1 m_1 e_1-\dots- \ell_{i-1} m_{i-1}
e_{i-1}$$
Recalling that $e_i=e_{i-2}-\ell_{i-1} e_{i-1}$, this simplifies to:
$$\boxed{
S^{j_c}\ell_i \int_{\loc B_{i-1}} \frac {k^{j_{\overline\mu}} 
(k-\ell_1 e_1-e_2)^{j_\mu} (k-e_1)^{j_\lambda} (\text{p.c.})^{j_c}}
{(1+\ell_i e_i)(1+e_{i-1})}\cap[\loc B_{i-1}]}$$\vskip 6pt

The total `local' contribution from the cusp $\lambda\cap\mu$ is
the sum of the first two boxes listed above, plus the sum of the last
box over $i=3,\dots,e$. Each of the terms can be evaluated
as a polynomial in $n$, the $m_i$'s, and the $\ell_i$'s, which we take
as indeterminates for a moment. With this understood, we let
$$Q_i(n,m_1,\ell_1,m_2,\ell_2,\dots,m_i,\ell_i)$$
be the contribution coming from the $i$@-th batch, and let
$$P_r(n,m_1,\ell_1,m_2,\ell_2,\dots,m_r,\ell_r)=\sum_{i=1}^r
Q_i(n,m_1,\ell_1,m_2,\ell_2,\dots,m_i,\ell_i)$$
The total local contribution from the point $\lambda\cap\mu$ is
then $P_e$, where the $\ell_i$'s are replaced with their values
prescribed by the ($e$@-line) Euclidean algorithm, that is:
$$P_e\left(n,m_1,\frac{n-m_2}{m_1},m_2,\frac{m_1-m_3}{m_2},\dots,
m_{e-1},\frac{m_{e-2}-m_e}{m_{e-1}},m_e,\frac{m_{e-1}}{m_e}\right)$$
Again, we can leave the $m_i$'s undetermined for a moment and treat
this as an expression in the variables $n,m=m_1,m_2,\dots,m_e$.
\proclaim{Lemma 3.4} As an expression in the variables
$n,m=m_1,m_2,\dots,m_e$,
$$P_e\left(n,m_1,\dots,\frac{m_{e-1}}{m_e}\right) = 
S^{j_c}n^{3-j_\lambda}\int_{\P^2} \frac
{k^{j_{\overline\mu}+1}(m+n k)^{j_c+j_\lambda-3}}{(1+k)^{3-j_\mu}}
\cap[\P^2]$$
\endproclaim
In particular, we are claiming that the above expression $P_e(\dots)$
does not depend on the intermediate multiplicities $m_2,\dots,m_e$,
and in fact it is independent on the number $e$ of lines taken by the
Euclidean algorithm on $m,n$. In a sense, all the needed information
is therefore contained {\it in nuce\/} in the contribution from the
{\it first\/} directed blow@-up!

\demo{Proof} The right@-hand@-side is $\dsize P_1\left(n,m,\frac
nm\right)$, that is the first box listed above but with $\ell_1=\frac
nm$. So the statement is correct for $e=1$; and the reader can check
it is correct for $e=2$. It can be proved for all
$e>2$ by induction. For this, it suffices to show that, for $e>2$,
the expression
$$P_e\left(n,m_1,\frac{n-m_2}{m_1},\dots,
m_{e-1},\frac{m_{e-2}-m_e}{m_{e-1}},m_e,\frac{m_{e-1}}{m_e}\right)$$
agrees with the expression
$$P_{e-1}\left(n,m_1,\frac{n-m_2}{m_1},\dots,
m_{e-1},\frac{m_{e-2}}{m_{e-1}}\right)$$
By definition:
$$\multline
P_e\left(n,m_1,\frac{n-m_2}{m_1},\dots,
m_{e-1},\frac{m_{e-2}-m_e}{m_{e-1}},m_e,\frac{m_{e-1}}{m_e}\right)\\
=P_{e-1}\left(n,m_1,\frac{n-m_2}{m_1},\dots,
m_{e-1},\frac{m_{e-2}-m_e}{m_{e-1}}\right)+Q_e\left(n,m_1,\dots,
m_e,\frac{m_{e-1}}{m_e}\right)
\endmultline$$
so it is enough to show that

(1) $\dsize Q_e\left(n,m_1,\dots, m_e,\frac{m_{e-1}}{m_e}\right)$ is a
polynomial in $m_e$, and vanishes for $m_e=0$; and

(2) $\dsize P_e\left(n,m_1,\dots, m_e,\frac{m_{e-1}}{m_e}\right)$ does
{\it not\/} depend on $m_e$.

Indeed, by (2) we may then assume $m_e=0$ in evaluating $P_e$; and by
(1) this will give
$$P_{e-1}\left(n,m_1,\frac{n-m_2}{m_1},\dots,
m_{e-1},\frac{m_{e-2}-0}{m_{e-1}}\right)+0$$
which is what is needed for the induction step.

Next, observe that the only summands in $P_e=\sum_{i=1}^e Q_i$ which
involve $m_e$ are $Q_{e-1}$ and $Q_e$; therefore in order to prove (2)
it is enough to prove

($2'$) $\dsize Q_{e-1}\left(n,m_1,\dots, m_{e-1},\frac{m_{e-2}-m_e}
{m_{e-1}}\right)+Q_e\left(n,m_1,\dots, m_e,\frac{m_{e-1}}{m_e}\right)$
does not depend on $m_e$.

Now we are only interested in the terms in $Q_{e-1}$, $Q_e$ involving
$m_e$, so we can neglect the homogeneous terms in $Q_i$ and absorb most
of the rest into a divisor class $D$. Then we are left with
$$\multline
\text{the term of codimension $j_c-2$ in} \\
\ell_{e-1} \frac {(m_{e-1}+D)^{j_c}}
{(1+\ell_{e-1} e_{e-1})(1+e_{e-2})}+
\ell_e \frac {(m_e+D- \ell_{e-1} m_{e-1} e_{e-1})^{j_c}}
{(1+\ell_e e_e)(1+e_{e-1})}
\endmultline$$
(where $\ell_{e-1}=\frac{m_{e-2}-m_e}{m_{e-1}}$ and
$\ell_e=\frac{m_{e-1}} {m_e}$). It is clear that the second summand is
a polynomial in $m_e$ and vanishes for $m_e=0$, as there are no terms
of codimension $j_c-2$ in this summand if $m_e=0$. Using \cite{A-F3}
and $e_e=e_{e-2}-\ell_{e-1} e_{e-1}$, the reader can check that the
sum equals the term in codimension $j_c-2$ in
$$\frac{m_{e-2}m_{e-1}(1+D)^{j_c}}{(1+m_{e-1}
e_{e-2})(1+m_{e-2}e_{e-1})}\quad,$$
so indeed it does not depend on $m_e$.\qed\enddemo

A simple substitution now computes the contribution due to the other
cusp, $\lambda\cap\overline\mu$. This amounts to replacing $m$ by
$n-m$, and of course reversing the roles of $\mu$ and $\overline\mu$.
Putting the two contributions together, combining with the
multiplicity data, and adding over $j_\mu$, etc. (cf.~Lemma~3.3) we
obtain the total `local' correction term:
$$\multline
7!\sum_{j_\mu+j_{\overline\mu}+j_\lambda+j_c=7}
n^{j_\mu+j_{\overline\mu}-4} \int_{\P^2} \left(\frac {k^{j_{\overline\mu+1}}
(m+n k)^{j_c+j_\lambda-3}}{(1+k)^{3-j_\mu}}\right.\\
\left. + \frac {k^{j_{\mu+1}}
(n-m+n k)^{j_c+j_\lambda-3}}{(1+k)^{3-j_{\overline\mu}}}\right)\,
\frac{q^{j_\mu}} {j_\mu!} \frac{\overline q^{j_{\overline\mu}}}
{j_{\overline\mu}!}\, \frac{r^{j_\lambda}}{j_\lambda!} \, \frac{(S
n)^{j_c}} {j_c!}\cap[\P^2]
\endmultline$$
where we are maintaining the convention that $q^3=\overline
q^3=r^3=0$.\vskip 6pt

This term can now be expanded with relative ease, yielding
$$\multline
n\,S^2(630\,m\,{q^2}\,{\overline q}\, {r^2} - 630\,m\,q\, {{{\overline
q}}^2}\,{r^2} + 630\,n\,q\, {{{\overline q}}^2}\,{r^2} +
420\,m\,n\,{q^2}\, {\overline q}\,r\,S - 420\,m\,n\,q\, {{{\overline
q}}^2}\,r\,S\\
+ 420\,{n^2}\,q\, {{{\overline q}}^2}\,r\,S -
210\,{m^2}\,{q^2}\,{r^2}\, S + 420\,m\,n\,{q^2}\, {r^2}\,S +
840\,{m^2}\,q\, {\overline q}\,{r^2}\,S - 840\,m\,n\,q\, {\overline
q}\,{r^2}\,S\\
+ 420\,{n^2}\,q\, {\overline q}\,{r^2}\,S - 210\,{m^2}\, {{{\overline
q}}^2}\, {r^2}\,S + 210\,{n^2}\, {{{\overline q}}^2}\, {r^2}\,S +
105\,m\,{n^2}\,{q^2}\, {\overline q}\,{S^2} - 105\,m\,{n^2}\,q\,
{{{\overline q}}^2}\,{S^2}\\
+ 105\,{n^3}\,q\, {{{\overline q}}^2}\,{S^2} -
105\,{m^2}\,n\,{q^2}\,r\, {S^2} + 210\,m\,{n^2}\, {q^2}\,r\,{S^2} +
420\,{m^2}\,n\,q\, {\overline q}\,r\,{S^2}\\
- 420\,m\,{n^2}\,q\, {\overline q}\,r\,{S^2} + 210\,{n^3}\,q\,
{\overline q}\,r\,{S^2} - 105\,{m^2}\,n\, {{{\overline q}}^2}\,r\,
{S^2} + 105\,{n^3}\, {{{\overline q}}^2}\,r\, {S^2} - 315\,{m^3}\,q\,
{r^2}\,{S^2}\\
+ 630\,{m^2}\,n\,q\,{r^2}\, {S^2} - 315\,m\,{n^2}\,q\, {r^2}\,{S^2} +
105\,{n^3}\,q\,{r^2}\, {S^2} + 315\,{m^3}\, {\overline
q}\,{r^2}\,{S^2}\\
- 315\,{m^2}\,n\, {\overline q}\,{r^2}\,{S^2} + 105\,{n^3}\,
{\overline q}\,{r^2}\,{S^2} - 21\,{m^2}\,{n^2}\,{q^2}\, {S^3} +
42\,m\,{n^3}\, {q^2}\,{S^3} + 84\,{m^2}\,{n^2}\,q\, {\overline
q}\,{S^3}\\
- 84\,m\,{n^3}\,q\, {\overline q}\,{S^3} + 42\,{n^4}\,q\, {\overline
q}\,{S^3} - 21\,{m^2}\,{n^2}\, {{{\overline q}}^2}\,{S^3} +
21\,{n^4}\, {{{\overline q}}^2}\,{S^3} - 126\,{m^3}\,n\,q\,r\, {S^3}\\
+ 252\,{m^2}\,{n^2}\, q\,r\,{S^3} - 126\,m\,{n^3}\,q\,r\,{S^3} +
42\,{n^4}\,q\,r\,{S^3} + 126\,{m^3}\,n\, {\overline q}\,r\,{S^3} -
126\,{m^2}\,{n^2}\, {\overline q}\,r\,{S^3}\\
+ 42\,{n^4}\,{\overline q}\, r\,{S^3} - 126\,{m^4}\,{r^2}\,{S^3} +
252\,{m^3}\,n\,{r^2}\, {S^3} - 126\,{m^2}\,{n^2}\, {r^2}\,{S^3} +
21\,{n^4}\,{r^2}\,{S^3}\\
 - 21\,{m^3}\,{n^2}\,q\,{S^4} + 42\,{m^2}\,{n^3}\,q\,{S^4} -
21\,m\,{n^4}\,q\,{S^4} + 7\,{n^5}\,q\,{S^4} + 21\,{m^3}\,{n^2}\,
{\overline q}\,{S^4}\\
- 21\,{m^2}\,{n^3}\, {\overline q}\,{S^4} + 7\,{n^5}\,{\overline q}\,
{S^4} - 42\,{m^4}\,n\,r\, {S^4} + 84\,{m^3}\,{n^2}\, r\,{S^4} -
42\,{m^2}\,{n^3}\,r\,{S^4}\\
+ 7\,{n^5}\,r\,{S^4} - 6\,{m^4}\,{n^2}\,{S^5} +
12\,{m^3}\,{n^3}\,{S^5} - 6\,{m^2}\,{n^4}\,{S^5} + {n^6}\,{S^5})
\endmultline$$
This expression is much more structured than it appears at first
sight. Using our convention ($q^3=\overline q^3=r^3=0$) we can rewrite
it as
$$\multline
(Sn+r+q+\overline q)^7-n^3 m^2 \overline m^2\left(\left(S+
\frac {r}n+\frac{q}m +\frac{\overline q}{\overline m}\right)^7+
2\left(S+\frac {r}n+\frac{q}m\right)^7 \right.\\
\left.+2\left(S+\frac {r}n+\frac{\overline q}{\overline m}\right)^7+
\left(S+\frac{r}n\right)^7-42 \left(S+\frac{r}n\right)^5
\left(\frac{q^2}{m^2}- \frac{q}{m} \frac{\overline q}{\overline m}+
\frac{\overline q^2}{\overline m^2}\right)\right)
\endmultline$$
(where $\overline m=n-m$). Subtracting from the B\'ezout term given in
the beginning of this subsection, we get the first expression listed in
Theorem~1.1:
$$\multline
n^3 m^2 \overline m^2\left(\left(S+ \frac {r}n+\frac{q}m
+\frac{\overline q}{\overline m}\right)^7+ 2\left(S+\frac
{r}n+\frac{q}m\right)^7 \right.\\
\left.+2\left(S+\frac {r}n+\frac{\overline q}{\overline m}\right)^7+
\left(S+\frac{r}n\right)^7-42 \left(S+\frac{r}n\right)^5
\left(\frac{q^2}{m^2}- \frac{q}{m} \frac{\overline q}{\overline m^2}+
\frac{\overline q^2}{\overline m^2}\right)\right)
\endmultline$$
(up to the multiplicative factor $1/(A n)$, cf.~Lemma 3.3).

This is the intersection product of the relevant divisors, after the
sequence of local blow@-ups is completed. The correction term for the
global blow@-ups is computed in the next subsection.

\subheading{\S3.2. Global blow@-ups}
Recall from \S2.2 and \S2.3
that after the local stages are completed, the base locus of the
rational map $c$ consists of several disjoint three@-dimensional
components, each isomorphic to $\P^2\times \P^1$, where the
$\P^1$@-factor represents the normalization of a curve of type
$(m,n)$, and of other five@-dimensional components due to the lines
$\lambda$, $\mu$, $\overline\mu$ of the basic triangle. 

We deal here with the three@-dimensional components. By Theorem~2.4,
two blow@-ups will resolve the indeterminacies of $c$ over such
components; we evaluate the relevant intersection products (as
in Lemma~3.3)
by subtracting from the result of \S3.1 a correction term due to these
components. The correction term will again be obtained by applying
Proposition~2.3,
for which we need to compile the usual information of multiplicity and
class of normal bundles. We will see that an interesting phenomenon
rules these `global' contributions: they are independent of $m$. 

First, we have to compute the Chern classes of the normal bundle to
each component $\P^2\times\P^1$; these will be 
$$\frac{c(T\Til\P^8)}{c(T \P^2\times \P^1)}$$
where $\Til \P^8$ is the blown@-up $\P^8$ at this stage, and we omit
the obvious pull@-back to $\P^2\times\P^1$. Now let $k,p$ denote the
pull@-back of the class of a line from the $\P^2$ factor and of a
point from the $\P^1$ factor; so
$$c(T \P^2\times \P^1)=(1+k)^3(1+2p)\quad.$$
As for $c(T\Til\P^8)$, this can be obtained by repeated application of
\cite{Fulton}, \S15.4: the reader will check that, after the sequence
of blow@-ups over $\lambda\cap\mu$, this pulls back to
$$\multline
(1+k+np)^9+(1+k+np)^8\cdot\\
((m_1+m_2)(k-2)p-3\ell_1m_1p-(\ell_2m_2+\ell_3 m_3+\dots)(k+1)p)
\endmultline$$
{From} the Euclidean algorithm, we see that
$$\ell_2 m_2+\ell_3
m_3+\dots=(m_1-m_3)+(m_2-m_4)+\dots+(m_{e-2}-m_e)+m_{e-1}$$
telescopes to $m_1+m_2-1$ since $m_e=gcd(m,n)=1$; and $\ell_1
m_1=n-m_2$. Hence the class simplifies to
$$\multline
(1+k+np)^9+(1+k+np)^8((m_1+m_2)(k-2)p-3(n-m_2)p-(m_1+m_2-1)(k+1)p)
\\
=c(T\P^8)+(1+k+np)^8((1-3n-3m)p+kp)
\endmultline$$
The second summand evaluates the change in the total class of the
tangent bundle due to the blow@-ups over the point $\lambda\cap\mu$.
Simply substituting $m\mapsto n-m$ evaluates the change due to the
other point, $\lambda\cap\overline\mu$:
$$(1+k+np)^8((1-3n-3(n-m))p+kp)$$
so that in total
$$\multline
c(T\Til\P^8)=(1+k+np)^9+(1+k+np)^8(((1-3n-3m)+(1-3n-3(n-m)))p+2kp)\\
=(1+k+np)^9+(1+k+np)^8((2-9n)p+2kp)
\endmultline$$
Expanding and applying $k^3=0$, $p^2=0$ (as we are pulling back to
$\P^2\times\P^1$) gives
$$c(T\Til\P^8)=1+9k+2p+36k^2+18 kp+72k^2p\quad:$$
remarkably, this expression is independent of $m,n$. In fact, the
normal bundle to one three@-dimensional component is computed by
$$\frac{c(T\Til\P^8)}{c(T \P^2\times \P^1)}= \frac{1+9k+2p+36k^2+18
kp+72k^2p}{(1+k)^3(1+2p)}=1+6k+15k^2=(1+k)^6$$
and is therefore particularly simple. This is the key ingredient in
evaluating the correction term due to the first global blow@-up.

As for the second, we use the analysis of the similar situation in
\cite{A-F1}
for smooth curves: indeed, the point of the proof of Theorem~2.4
is that after the local blow@-ups and the first global blow@-up,
the geometry over one of these components is entirely analogous to the
geometry for smooth curves, over non@-flex points.
In particular, we know that the center of the second global blow@-up
is the union of a $\P^1$@-bundle over each three@-dimensional
component of the center of the first blow@-up, with $c_1(\Cal
O(-1))=f=$restriction of the exceptional divisor, and standard
computations
will show that, after push@-forward to the underlying $\P^2\times\P^1$,
$$f^4\mapsto 0,\, f^3\mapsto -6k^2,\, f^2\mapsto -3k,\, f\mapsto -1\quad.$$
The normal bundle has class
$$(1+f)(1+k-f)^3\quad.$$

The last ingredient necessary to perform the computation of the global
contribution is the restriction to the centers of the classes of the
proper transforms of the divisors. To find these, keep in mind that
$e_i$ restricts to $m_i p$; then

---the proper transform of a general hyperplane $\Til H$ restricts to
$(k+np)$;

---the proper transform of $H_\lambda$ restricts to
$k+np-e_1-\overline e_1$ (where $e_1$, $\overline e_1$ are the first
exceptional divisors  at $\lambda\cap\mu$, $\lambda\cap\overline\mu$,
resp.) This equals $k+np-mp-(n-m)p=k$;

---the proper transform of $H_\mu$ restricts to $k+np-\ell_1e_1-e_2=
k+np-(\ell_1 m_1+m_2)=k+np-np=k$; similarly, the proper transform of
$H_{\overline\mu}$ must restrict to $k$;

the proper transform of $X$ restricts to 
$$Sn(k+n p)-m_1 S \ell_1 e_1-m_2 S \ell_2 e_2 - \cdots$$
after the sequence over the point $\lambda\cap\mu$, that is to
$$\align
Sn(k+np) &-S(\ell_1 m_1^2-\ell_2 m_2^2-\cdots)p\\
&=Sn(k+np)-S(m(n-m_2)+m_2(m-m_3)+\cdots)p\\
&=Sn(k+np)-Sm n p\quad.
\endalign$$
The second summand is the change due to the sequence of blow@-ups
over $\lambda\cap\mu$; to obtain the change due to $\lambda\cap
\overline\mu$, just substitute $m\mapsto n-m$. The conclusion is that

---the proper transform of $X$ at the first global blow@-up restricts to 
$$Sn(k+np)-Smnp-S(n-m)np=Snk\quad.$$
Also, note that $X$ has multiplicity $s_i$ along the $i$@-th
component, 0 along all others.

Note that none of the classes depend on $m$: as we claimed above, {\it
the global contributions do not depend on $m$.\/}

Putting the above together and using Proposition~2.3,
we see that the global contributions are obtained by evaluating
$$\boxed{\sum_i \int \frac{(k+np) k^{7-j_c}
(s_i+ Snk)^{j_c}} {(1+k)^6}}$$
and
$$\boxed{\sum_i \int \frac{(k+np) k^{7-j_c} (s_i+ Snk-s_if)^{j_c}}
{(1+f)(1+k-f)^3}}\quad.$$
Adding these two terms and inserting in Lemma~3.3
gives a relatively simple expression:
$$
n\,\left( 84\,{{\left( Sn + r + q + \overline q \right)}^2}\,{\sum
s_i^5} - 252\,\left(S n + r + q + \overline q \right)\,{\sum s_i^6} +
192\,{\sum s_i^7} \right)
$$
which reproduces the one given in \S1, again up to the multiplicative
factor $1/(An)$.

\head \S4. End of the computation, and variations\endhead

The careful reader knows that we are not quite done, since after the
pair of global blow@-ups we are still left with base loci
corresponding to the lines of the basic triangle (cf.~\S2.3).
 
Here we reap the benefit of having shown that we only need to compute
the relevant intersection products for $j_\mu$, $j_{\overline\mu}$,
$j_\lambda\le 2$ (Lemma~3.3). Indeed, with
these constraints the corrections due to these base loci are {\it
zero.\/} For example, consider the correction term coming from the
$\P^2$ of matrices whose image is the point $\mu\cap\overline\mu$:
denoting by $k$ the class of a line in $\P^2$ and applying once more
the formula in Proposition~2.3,
this term is evaluated as
$$\int_{\P^2} \frac{k(q+qk)^{j_\mu} (\overline q+\overline q
k)^{j_{\overline\mu}} (rk)^{j_\lambda} (Snk)^{j_c}} {(1+k)^6}\quad:$$
and since $1+j_c+j_\lambda=8-(j_\mu+j_{\overline\mu})\ge 4$, this term
is automatically~0.
 
The same discussion applies to the remaining three 5@-dimensional
base loci; we leave the details to the reader. Theorem~1.1 then
follows,
since this shows that the expressions obtained in \S3.1 and \S3.2
combine to give the intersection product in Lemma~3.3.

\subheading{\S4.1. Quadritangent conics}
There is (in characteristic zero) only one 
other class of curves $C$ whose components are not all lines and whose
orbit is small: $C$ consists of 2 or more conics from a pencil through
a conic and a double tangent line; it may also contain that tangent
line. The multiplicities of components are arbitrary. For this class
of curves, the stabilizer is 1@-dimensional; its identity component
is the additive group $\G_a$.
$$\epsffile{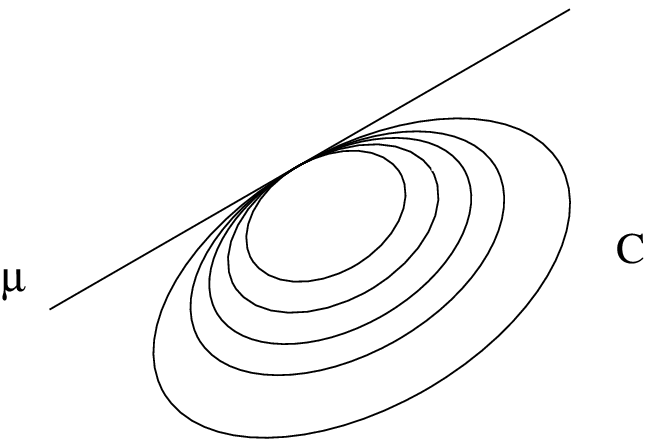}$$

These curves are not of the type considered in previous sections.  A
variety $\Til V$ dominating the orbit closure of such a curve can
however be constructed by a strategy very similar to the one followed
in \S2: again, the sequence of blow@-ups giving an embedded resolution
of the curve can be mirrored to produce a variety $\loc{\Til V}$, from
which a variety $\Til V$ is produced by the technique of Theorem~2.4.
Intersection@-theoretic computations similar to those in \S3 allow us
then to compute the degree of the orbit closure of these curves.

It is perhaps a little surprising that the formula given in
Theorem~1.1
turns out to be correct for this case as well: taking $n=2$,
$m=\overline m=1$, and $\overline q=r = 0$ computes the degree of the
orbit closure of a union of quadritangent conics appearing with
multiplicity $s_i$, together with the tangent line $\mu$ at the
point of contact, taken with multiplicity $q$. In other words, the
polynomial $Q$ for such a curve is the same as the polynomial for a
union of `bitangent conics'. 
Again, $A$ equals the number of components of the stabilizer of $C$; 
the analogue of Lemma~3.1
is the fact that for $C$ given by
$y^q \prod_i (x^2+yz+\alpha_i y^2)^{s_i}=0$, the number $A$
equals twice the maximum order of an automorphism 
$\alpha\mapsto u\alpha+v$ preserving the $S$@-tuple given by
$\prod(\alpha-\alpha_i)^{s_i}=0$.
Hence $A=2$ in most cases, and it is bounded by twice the number of conics
in $C$; for two conics with equal multiplicities, $A=4$.

With notations as above, the polynomial $Q$ is in this case
$$
24 S^7+84S^6q+84S^5q^2-84(2S+q)^2\sum s_i^5
+252(2S+q)\sum s_i^6-192\sum s_i^7
$$
For example, there are 504(=2016/4) pairs of quadritangent conics
through 7 general points.

Expressions for the degree of loci corresponding to curves with fixed
tangent line, or tangent line constrained to contain a given point,
can be obtained by differentiating this expression with respect to $q$
(cf.~\S4.3).

\subheading{\S4.2. Predegree polynomials} Simple adjustments in the
computations described in this paper allow us to compute the degrees
of suitable subsets of the orbit closure, obtained by imposing general
linear conditions on the matrices used to act on $C$. In a sense this
note deals precisely with one such computation: we computed the degree
of $\overline{\Cal O}_C$ by imposing a general linear condition on
$\P^8$ (and arguing that this would intersect the fibers over a point
of $\Cal O_C$ in $A n$ points, cf.~Proposition~3.2).
We call the {\it `predegree'\/} of $\overline{\Cal O}_C$ the product
of the degree of $\Cal O_C$ with the degree in $\P^8$ of the closure
of the stabilizer of $C$. Arguing as in the discussion leading to
Proposition~3.2 (and using the same notations),
we see that if $\Cal O_C$ has dimension $k$ then
$$\text{(predegree of $\Cal O_C$)}=\int \Til H^{8-k} \Til W^k\quad,$$
where the intersection product is taken in any variety resolving the
indeterminacies of the relevant rational map. (Note: this notion of
predegree agrees with the one used in \cite{A-F1}, where $k=8$.)

We find it in fact useful to introduce an {\it `adjusted predegree
polynomial'} defined by
$$\sum_{j\ge 0} \left(\int \Til H^{8-j} \Til W^j\right)
\frac{t^j}{j!}\quad;$$
if $\dim\Cal O_C=k$, then the coefficient of $t^k/k!$ in the adjusted
predegree polynomial gives the predegree of $\Cal O_C$, while for
$j>k$ the coefficient of $t^j/j!$ is~0. 

The introduction of denominators reflects some extra structure of
these polynomials, which we will not discuss here. As an example, note
the factorization of the polynomial for the degenerate case of a curve
supported on the basic triangle, indicated below: such factorizations
would not occur `without denominators'.

Suitable variations of the computations in \S3 yield the whole
adjusted predegree polynomial for the curves considered in this paper; the
result is as follows. To obtain the polynomial for a curve with data
$m,n$, $\overline m=n-m$, $s_i$, $S=\sum s_i$, $r,q,\overline q$, as
above: {\it imposing $q^3=\overline q^3=r^3=0$ in all computations,\/}
truncate to $t^7$ the expansion of 
$$e^{(S n+r+q+\overline q)t}\quad,$$
(the `B\'ezout' term), then subtract a global contribution
(independent of $m$){\eightpoint
$$\multline
12\,n\,{\sum s_i^5}\frac{t^5}{5!} + n\,\left( 48\,(S n+r+q+\overline
q)\,{\sum s_i^5} - 72\,{\sum s_i^6} \right)\frac{t^6}{6!}\\
+ n\,\left( 84\,{(S n+r+q+\overline q)^2}\,{\sum s_i^5} - 252\,(S n+
r+ q+ \overline q) \,{\sum s_i^6} + 192\,{\sum s_i^7} \right)\frac{t^7}{7!}
\endmultline$$}
and a local contribution given in degree~6 by
$$n^3 \left(m^3\left(S+\frac rn+\frac qm\right)^6+
\overline m^3 \left(S+\frac rn+\frac {\overline q}{\overline m}\right)^6
\right)\frac{t^6}{6!}$$
and in degree 7 by:
$$\multline
\left((Sn+r+q+\overline q)^7-n^3 m^2 \overline m^2\left(\left(S+
\frac {r}n+\frac{q}m +\frac{\overline q}{\overline m}\right)^7+
2\left(S+\frac {r}n+\frac{q}m\right)^7 \right.\right.\\
\left.\left.+2\left(S+\frac {r}n+\frac{\overline q}{\overline m}\right)^7+
\left(S+\frac{r}n\right)^7-42 \left(S+\frac{r}n\right)^5
\left(\frac{q^2}{m^2}- \frac{q}{m} \frac{\overline q}{\overline m}+
\frac{\overline q^2}{\overline m^2}\right)\right)\right) \frac {t^7}{7!}
\endmultline$$

The advantage of looking at the whole polynomial is that it carries
degree information for all orbits, regardless of their dimension (while
Theorem~1.1
does assume that the orbit of the curve under exam
has dimension~7). For example, setting $m=1$, $n=2$, $S=s_1=1$, and
$r=q=\overline q=0$ gives a polynomial
$$1+2t+\frac{4t^2}2+\frac{8t^3}{3!}+\frac{16 t^4}{4!} + \frac{8
t^5}{5!}\quad; $$
as the orbit closure of a conic is clearly the whole of $\P^5$, this
correctly detects that the stabilizer of a conic is a threefold of
degree~8. In fact, the second Veronese embedding of the $\P^3$ of
$2\times 2$ matrices in the
$\P^9$ of space quadrics projects isomorphically to this threefold in
$\P^8$; with suitable identifications, the center of the projection is
the determinant quadric.

For another example, take all $s_i=0$: that is, consider a curve
consisting solely of lines supported on the sides $\mu$,
$\overline\mu$, $\lambda$ of the basic triangle, with multiplicities
$q$, $\overline q$, $r$. This yields a degree@-6 polynomial, which in
fact factors
$$\left(1+q t+\frac{q^2 t^2}2\right) \left(1+\overline q t+\frac{
\overline q^2 t^2}2\right) \left(1+r t+\frac{r^2 t^2}2\right)\quad:$$
the orbit has dimension~6 and predegree $90 q^2\overline q^2 r^2$, as
the reader could check independently by observing that this orbit
closure can also be realized as the image of the evident map $\P^2
\times \P^2\times \P^2 @>>> \P^N$ ($N=d(d+3)/2$ for $d=q+\overline
q+r$).
The degree of the closure of the stabilizer depends on the
multiplicities of the lines: it is $3!$ if $q=\overline q=r$, $2$ if
exactly two multiplicities agree, and $1$ if the multiplicities are
distinct. 

All such computations are very particular cases of the general
expression for the adjusted predegree polynomial given above. This
covers then almost all small orbits of plane curves, with few
exceptions such as curves consisting of a star of lines through a
point $p$, union a line not containing $p$. Note that these are curves
`of type $(1,1)$' according to the terminology used in this paper; but
the blow@-up construction of \S2 assumes that the only linear
components of the curve are the lines of the basic triangle, so the
construction fails in this case. We will consider these curves
in \cite{A-F4}.

\subheading{\S4.3. Curves with constraints} The
polynomial given in Theorem~1.1
also contains enumerative information on the subsets of the orbits
parametrizing curves with specified constraints on the lines of the
basic triangle (such as: containing a given point). Let $\overline{\Cal
O}_C (j_\mu,j_{\overline\mu},j_\lambda)$ be the closure of the set of
translations $C\circ\varphi$ such that $\mu\circ\varphi$ contains
$j_\mu$ given points, $\overline\mu\circ\varphi$ contains
$j_{\overline\mu}$ given points, and $\lambda\circ\varphi$ contains
$j_\lambda$ given points (all choices of the points being general).

\proclaim{Proposition 4.1} Let $Q(n,m,s_i,r,q,\overline q)/A$ be the
polynomial giving the degree of $\overline{\Cal O}_C$, as in
Theorem~1.1.
Then the degree of $\overline{\Cal O}_C
(j_\mu,j_{\overline\mu},j_\lambda)$ is
$$\frac{(7-j_\mu- j_{\overline\mu} -j_\lambda)!}{7!A}
\frac{\partial^{j_\mu}}{\partial q^{j_\mu}}
\frac{\partial^{j_{\overline\mu}}}{\partial \overline q^{j_{\overline\mu}}}
\frac{\partial^{j_\lambda}}{\partial r^{j_\lambda}}
Q(n,m,s_i,r,q,\overline q)\quad.$$
\endproclaim
\demo{Proof} Arguing as in Proposition~3.2,
$$A\,n\,\deg\overline{\Cal O}_C
(j_\mu,j_{\overline\mu},j_\lambda)=\int \Til H\cdot \Til W^{(7-
j_\mu-j_{\overline\mu}-j_\lambda)}\cdot \Til H_\mu^{j_\mu} \cdot\Til
H_{\overline\mu}^{j_{\overline\mu}}\cdot \Til H_{\lambda}^{j_\lambda}
$$
where $\Til W=q\Til H_\mu+\overline q\Til H_{\overline\mu} +r\Til
H_\lambda +\Til X$ is the class of a point@-condition in a variety
resolving the rational map corresponding to $C$. Now a direct
computation shows that this equals
$$\frac{(7-j_\mu- j_{\overline\mu} -j_\lambda)!}{7!}
\frac{\partial^{j_\mu}}{\partial q^{j_\mu}}
\frac{\partial^{j_{\overline\mu}}}{\partial \overline q^{j_{\overline\mu}}}
\frac{\partial^{j_\lambda}}{\partial r^{j_\lambda}}\int
\Til H\cdot (q\Til H_\mu+\overline q\Til H_{\overline\mu} +r\Til 
H_\lambda +\Til X)^7\quad,$$
which gives the statement.\qed\enddemo

This result has a clear enumerative meaning for example when
$q=\overline q=r=0$, that is when $C$ does not contain the lines in
the triangle, and when all $s_i=1$ (so that $S=$ the number of
components of $C$). Then the 
formula of Proposition~4.1
gives  the number of $\PGL(3)$-translations of
$C$ which satisfy the given constraints on the lines of the basic
triangle, and contain the appropriate number ($=7-j_\mu-j_{\overline
\mu}-j_\lambda$) of general points.

{From} the formula in Theorem~1.1
and Proposition~4.1
it is easy to obtain closed formulas for these numbers. For example,
the number of unconstrained curves of fixed type (that is, $S$ given
distinct components from the pencil of type $(m,n)$@-curves)
through $7$ points is
$$\frac{6 S}A \left(m^2 \overline m^2 n^2 S^6-14 n^2 S^2+42
n S-32\right)\quad;$$
for the simplest example, take $n=3$, $m=2$ and $S=1$: there are 24
cuspidal plane cubics through 7 general points.
The number of curves such that the line $\lambda$ contains a given
point is
$$\frac {6 S}A \left(m^2 \overline m^2 n S^5-4n S+6\right)\quad;$$
for example, there are 36 cuspidal cubics through 6 general points,
and such that the line connecting the flex and the cusp contains a
given point.

Analogously, the number of curves such that $\lambda$ contains two
given points
is, according to Proposition~4.1:
$$\frac{2 S}A \left(3 m^2 \overline m^2 S^4-2\right)\quad;$$
hence, there are 20 cuspidal cubics through 5 points, with fixed line
through cusp and flex.

The reader who so wishes will have no difficulty deriving the other 24
closed formulas for assortments of conditions on the lines of the basic
triangle. For cuspidal cubics, the 27 numbers so obtained reproduce
results in \cite{M-X} (where hundreds more are computed). Here they
are, where the number at the $(i,j,k)$ spot denotes the number of
curves with $i$ points through $\mu$, $j$ points through
$\overline\mu$, and $k$ points through $\lambda$:
$$
x^3=y^2\qquad:\qquad
\spreaddiagramrows{-2pc} \spreaddiagramcolumns{-2pc}
\diagram
 & & & & 4 \xdotted '[1,-2] \xdotted '[0,3] \xdotted '[1,0]'[2,0]'[3,0] & &
& 2 \xdotted '[1,-2] \xdotted '[0,3] \xdotted '[1,0]'[3,0] & & & 1 \xdotted
'[1,-2] \xdotted '[3,0] \\
 & & 12 \xdotted '[1,-2] \xdotted '[0,3] \xdotted '[1,0]'[3,0] & & & 4 \xdotted
'[1,-2] \xdotted '[0,3] \xdotted '[1,0]'[3,0] & & & 1 \xdotted '[1,-2]
\xdotted '[3,0] \\
20 \xdotted '[0,3] \xdotted '[3,0] & & & 6 \xdotted '[0,3] \xdotted '[3,0] & &
& 1 \xdotted '[3,0] \\
 & & & & 12 \xdotted '[1,-2] \xdotted '[0,1]'[0,2]'[0,3] \xdotted
'[1,0]'[2,0]'[3,0] & & & 6 \xdotted '[1,-2] \xdotted '[0,1]'[0,3] \xdotted
'[1,0]'[3,0] & & & 3 \xdotted '[1,-2] \xdotted '[3,0] \\
 & & 32 \xdotted '[1,-2] \xdotted '[0,1]'[0,3] \xdotted '[1,0]'[3,0] & & & 12
\xdotted '[1,-2] \xdotted '[0,1]'[0,3] \xdotted '[1,0]'[3,0] & & & 3 \xdotted
'[1,-2] \xdotted '[3,0] \\
36 \xdotted '[0,3] \xdotted '[3,0] & & & 14 \xdotted '[0,3] \xdotted '[3,0] & &
& 3 \xdotted '[3,0] \\
 & & & & 32 \xdotted '[1,-2] \xdotted '[0,1]'[0,2]'[0,3] & & & 18 \xdotted
'[1,-2] \xdotted '[0,1]'[0,3] & & & 9 \xdotted '[1,-2] \\
 & & 72 \xdotted '[1,-2] \xdotted '[0,1]'[0,3] & & & 32 \xdotted '[1,-2]
\xdotted '[0,1]'[0,3] & & & 9 \xdotted '[1,-2]\\
24 \xdotted '[0,3] & & & 18 \xdotted '[0,3] & & & 5\\
\enddiagram$$
For higher degree curves, the results are, to our knowledge, new.
Here are the numbers for the curve of (randomly chosen) type $(7,4)$:
$$
x^7=y^4\qquad:\qquad
\spreaddiagramrows{-2pc} \spreaddiagramcolumns{-2pc}
\diagram
 & & & & 16 \xdotted '[1,-2] \xdotted '[0,3] \xdotted '[1,0]'[2,0]'[3,0] & &
& 4 \xdotted '[1,-2] \xdotted '[0,3] \xdotted '[1,0]'[3,0] & & & {\,\,1\,\,} \xdotted
'[1,-2] \xdotted '[3,0] \\
 & & 144 \xdotted '[1,-2] \xdotted '[0,3] \xdotted '[1,0]'[3,0] & & &
24 \xdotted
'[1,-2] \xdotted '[0,3] \xdotted '[1,0]'[3,0] & & & 3 \xdotted '[1,-2]
\xdotted '[3,0] \\
860 \xdotted '[0,3] \xdotted '[3,0] & & & 108 \xdotted '[0,3] \xdotted '[3,0] & &
& 9 \xdotted '[3,0] \\
 & & & & 112 \xdotted '[1,-2] \xdotted '[0,1]'[0,2]'[0,3] \xdotted
'[1,0]'[2,0]'[3,0] & & & 28 \xdotted '[1,-2] \xdotted '[0,1]'[0,3] \xdotted
'[1,0]'[3,0] & & & {\,\,7\,\,} \xdotted '[1,-2] \xdotted '[3,0] \\
 & & 1004 \xdotted '[1,-2] \xdotted '[0,1]'[0,3] \xdotted '[1,0]'[3,0] & & & 168
\xdotted '[1,-2] \xdotted '[0,1]'[0,3] \xdotted '[1,0]'[3,0] & & & 21 \xdotted
'[1,-2] \xdotted '[3,0] \\
5916 \xdotted '[0,3] \xdotted '[3,0] & & & 752 \xdotted '[0,3] \xdotted '[3,0] & &
& 63 \xdotted '[3,0] \\
 & & & & 780 \xdotted '[1,-2] \xdotted '[0,1]'[0,2]'[0,3] & & & 196 \xdotted
'[1,-2] \xdotted '[0,1]'[0,3] & & & {\,\,49\,\,} \xdotted '[1,-2] \\
 & & 6924 \xdotted '[1,-2] \xdotted '[0,1]'[0,3] & & & 1172 \xdotted '[1,-2]
\xdotted '[0,1]'[0,3] & & & 147 \xdotted '[1,-2]\\
39792 \xdotted '[0,3] & & & 5160 \xdotted '[0,3] & & & 437\\
\enddiagram$$

\leftheadtext{References}
\rightheadtext{References}

\enddocument